\documentclass[12pt]{article}
\usepackage{enumitem}
\usepackage{fullpage}
\usepackage[english]{babel}
\usepackage[utf8]{inputenc}
\usepackage[authoryear,comma]{natbib}
\usepackage{amsmath,multibib}
\usepackage{bbm}
\usepackage{tabularx}

\usepackage{amssymb,amsthm,mathabx,mathtools}
\usepackage{mathrsfs}
\usepackage{graphicx}
\usepackage{url}
\usepackage{caption,subcaption}
\usepackage{booktabs}
\usepackage{multirow}
\usepackage{caption}
\usepackage{epsfig,float}
\usepackage{url} 
\usepackage[usenames,dvipsnames]{color} 
\usepackage{listings}
\usepackage{mathrsfs,bm} 

\usepackage{tikz}
\usepackage{graphicx}

\theoremstyle{plain}
\newtheorem{thm}{Theorem}

\newcommand{\bthm}{\begin{thm}}
\newcommand{\ethm}{\end{thm}}
\newcommand{\bpf}{\begin{proof}}
\newcommand{\epf}{\end{proof}}
\theoremstyle{definition}
\newtheorem{defn}{Definition}
\newtheorem{example}{Example}
\newtheorem{rem}{Remark}
\usepackage[margin=.9in,a4paper]{geometry}
\usepackage[compact]{titlesec}
    \titlespacing{\section}{0pt}{1ex}{1ex}
    \titlespacing{\subsection}{0pt}{1ex}{0ex}
    \titlespacing{\subsubsection}{0.5pt}{0.55ex}{0ex}
\numberwithin{equation}{section}
\usepackage{deep-macro}
\usepackage{epigraph}
\usepackage{relsize}
\usepackage{multicol}

 \newcommand\scalemath[2]{\scalebox{#1}{\mbox{\ensuremath{\displaystyle #2}}}}

\title{MaxEnt-Copula-Oct2020}
\author{mesuvodeep }
\date{September 2020}

\usepackage{deep-macro}
\usepackage{cancel}

\begin{document}
\begin{center} 
{\bf 
{\Large A Maximum Entropy Copula Model for Mixed Data:\\[.2em]
Representation, Estimation, and Applications}\\[1em]
Deep Mukhopadhyay}\\[.1em]
\texttt{deep@unitedstatalgo.com}\\[.4em]
{\bf Final Version: August, 2022}\\[.12em]
\end{center}
\begin{abstract} 
A new nonparametric model of maximum-entropy (MaxEnt) copula density function is proposed, which offers the following advantages: (i) it is
valid for \textit{mixed} random vector. By `mixed' we mean the method works for any combination of discrete or continuous variables in a \textit{fully automated} manner; (ii) it yields a \textit{bonafide} density estimate with \textit{intepretable} parameters. By `bonafide' we mean the estimate guarantees to be a non-negative function, integrates to $1$; and (iii) it plays a \textit{unifying role} in our understanding of a large class of statistical methods for mixed $(X,Y)$. Our approach utilizes modern machinery of nonparametric statistics to represent and approximate log-copula density function via LP-Fourier transform. Several real-data examples are also provided to explore the key theoretical and practical implications of the theory.
\end{abstract} 
\vspace{-.4em}
\linespread{1.27}
\renewcommand{\baselinestretch}{1.34}
\setlength{\parskip}{1.5ex}

\section{Copula Statistical Learning}
Copulas are the `bridge’ between the univariate and the multivariate statistics world, with applications in a wide variety of science and engineering fields---from economics to finance to marketing to healthcare. Because of the ubiquity of copula in empirical research, it is becoming necessary to develop a general theory that can unify and simplify the copula learning process. In this paper, we present a new class \textit{specially-designed} nonparametric maximum-entropy (MaxEnt) copula model
that offers the following advantages: First, it yields a \textit{bonafide} (smooth, non-negative, and integrates to $1$) copula density estimate with \textit{interpretable} parameters that provide insights into the nature of the dependence between the random variables $(X, Y)$. Secondly, the method is data-type agnostic---which is to say that it \textit{automatically} (self) adapts to mixed-data types (any combination of discrete, continuous, or even categorical). Thirdly, and most notably, our copula-theoretic framework subsumes and \textit{unifies a wide range of statistical learning methods} using a common mathematical notation---unlocking deep, surprising connections and insights, which were previously unknown.  In the development of our theory and algorithms, the LP-Fourier method of copula modeling (which was initiated by \cite{D20copula}), plays an indispensable role.

\section{Self-Adaptive Nonparametric Models}
\label{sec:theory}
We introduce two new classes of  maximum-entropy (MaxEnt) copula density models. But before diving into technical details, it will be instructive to review some basic definitions and concepts related to copula. 
\subsection{Background Concepts and Notation} \label{sec:back}
\vskip.25em
{\bf Sklar’s Copula Representation Theory} \citep{sklar1959}. The joint cumulative distribution function (cdf) of any pair of random variables $(X,Y)$
\[F_{X,Y}(x,y)=\Pr(X\le x, Y\le y),~~\mathrm{for}~(x,y) \in \cR^2\]
can be decomposed as a function of the marginal cdfs $F_X$ and $F_Y$
\beq \label{eq:Copcdf}
F_{X,Y}(x,y)\,=\,\Cop_{X,Y}\big( F_X(x), F_Y(y)\big), ~~\mathrm{for}~(x,y) \in \cR^2
\eeq
where $\Cop_{X,Y}$ denotes a copula distribution function with uniform marginals. To set the stage, we start with the continuous marginals case, which will be generalized later to allow mixed-$(X,Y)$. Taking derivative of Eq. \eqref{eq:Copcdf}, we get
\beq \label{eq:coppdf}
f_{X,Y}(x,y)\,=\,f_X(x) f_Y(y)\cop_{X,Y}\big( F_X(x), F_Y(y)\big), ~~\mathrm{for}~(x,y) \in \cR^2
\eeq
which decouples the joint density into the marginals and the copula. One can rewrite Eq. \eqref{eq:coppdf} to represent copula as a ``normalized'' joint density function 
\beq \label{eq:depf}
\cop_{X,Y}\big( F_X(x), F_Y(y)\big)\,:=\,{\rm dep}_{X,Y}(x,y)\,=\dfrac{f_{X,Y}(x,y)}{f_X(x) f_Y(y)},~~~~~~~~~~~
\eeq 
which is also known as the dependence function, pioneered by \cite{Hoeff40}. To make \eqref{eq:depf} a proper density function (i.e., one that integrates to one) we perform quantile transformation by substituting
$F_X(x)=u ~\text{and}~ F_Y(y)=v$:
\[
\cop_{X,Y}(u,v)={\rm dep}_{X,Y}( F_X^{-1}(u),  F_Y^{-1}(v))=\dfrac{f_{X,Y}\big( F_X^{-1}(u), F_Y^{-1}(v)  \big)}{f_X( F_X^{-1}(u) ) f_Y(F_Y^{-1}(v) )}, ~~(u,v) \in [0,1]^2.
\]
We are now ready to extend this copula density concept to the mixed $(X,Y)$ case.
\vskip.5em
{\bf Pre-Copula: Conditional Comparison Density} \citep{D13a}. Before we introduce the generalized copula density, we need to introduce a new concept---conditional comparison density (CCD). For a continuous $X$, CCD is defined as:
\beq  \label{eq:du}
d(u;X,X|Y=y)\,=\,\dfrac{f_{X|Y}\big( F_X^{-1}(u) \mid y  \big)}{f_X\big( F_X^{-1}(u)\big)},~~0<u<1 \eeq
For $Y$ discrete, we represent it using probability mass function (pmf):
\beq \label{eq:dv}
d(v;Y,Y|X=x)\,=\,\dfrac{p_{Y|X}\big( Q_Y(v)|x  \big)}{p_Y\big( Q_Y(v)\big)}\,=\,\dfrac{\Pr(Y=Q_Y(v)|X=x)}{\Pr(Y=Q_Y(v))},~~0<v<1. \eeq
where $Q_Y(v)$ is the quantile function of $Y$. It is easy to see that the CCDs \eqref{eq:du} and \eqref{eq:dv} are proper densities in the sense that
\[\int_0^1  d(u;X,X|Y=y)\dd u\,=\,\int_0^1 d(v;Y,Y|X=x) \dd v\,=\,1.\]
\vskip.12em
{\bf Generalized Copula Representation Theory} \citep{D20copula}. For the mixed case, when $Y$ is discrete and $X$ is continuous the joint density of \eqref{eq:depf} is defined by either side of the following identity:
\[
\Pr(Y\mid X=x) f_X(x)\,=\,f_{X|Y}(x|y) \Pr(Y=y).
\]
This can be rewritten as the ratios of conditionals and their respective marginals:
\beq \label{eq:pbr}
\text{Pre-Bayes' Rule}:~~~~\dfrac{\Pr(Y=y\mid X=x)}{\Pr(Y=y)}\,=\,\dfrac{f_{X|Y}(x|y)}{f_X(x)}~~~~~~~~~~
\eeq
This formula \eqref{eq:pbr} can be interpreted as the \textit{slices} of the \textit{mixed} copula density, since 
\beq \label{eq:copC}
\cop_{X,Y}\big( F_X(x), F_Y(y) \big)\,=\,\dfrac{f_{X|Y}(x|y)}{f_X(x)} \,=\,\dfrac{\Pr(Y=y\mid X=x)}{\Pr(Y=y)}.~~~~~
\eeq
Substituting $F_X(x)=u$ and $F_Y(y)=v$, we get the following definition of the  generalized copula density in terms of conditional comparison density (CCD): 
\beq \label{eq:copG}
\cop_{X,Y}(u, v)\,=\,d\big(u;X,X|Y= Q(v; Y)\big)\,=\,d\big(v;Y, Y|X=Q(u;X)\big),~~0<u,v<1.\eeq
Bayes' theorem ensures the equality of two CCDs, with copula being the common value.  Equipped with this fundamentals, we now develop the nonparametric theory of MaxEnt copula modeling.


\subsection{Log-bilinear Model}
An exponential Fourier series representation of copula density function is given. The reasons for entertaining an exponential model for copula
is motivated from two different perspectives.

\vskip.35em
{\bf The Problem of Unboundedness}. 
One peculiar aspect of copula density function is that it can be unbounded at the corners of the unit square. In fact, many common parametric copula families---Gaussian, Clayton, Gumbel, etc.---tend to infinity at the boundaries. So naturally the question arises: How to develop suitable approximation methods that can accommodate a  broader class of copula density shapes, including the unbounded ones? The first key insight: logarithm of the copula density function is far more convenient to approximate (due to its well-behaved nature) than the original density itself. We thus express the logarithm of copula density $\log \cop_{X,Y}$ in the Fourier series---instead of doing canonical $L_2$ approximation, which expands $\cop_{X,Y}$ directly in an orthogonal series \citep{D20copula}. Accordingly, for densities with rapidly changing tails, `log-Fourier' method leads to an improved estimate that is less wiggly and more parsimonious than the $L_2$-orthogonal series model. 
In addition, the resulting exponential form guarantees the non-negativity of the estimated density function. 
\vskip.1em

\textit{Choice of Orthonormal Basis}. To expand log-copula density function, we choose the LP-family of polynomials (see Appendix \ref{app:LPb}), which are especially suited to approximate functions of mixed random variables. In particular, we approximate $\log \cop_{X,Y}$ by expanding it in the tensor-product of LP-bases $\{S_j \otimes S_k\}$, which are orthonormal with respect to the empirical-product measure  $\{\wtF_X \otimes \wtF_Y\}$. LP-bases' appeal lies in its ability to approximate the quirky shapes of mixed-copula functions in a completely automated way; see Fig. \ref{fig:copFig1}. Consequently, it provides a unified way to develop nonparametric smoothing algorithms that simultaneously hold for mixed data types.
\vskip.3em
\begin{defn}
The exponential copula model admits the following LP-expansion
\beq \label{exp:model}
\cop_{\teb}(u,v;X,Y)~=~\dfrac{1}{Z_{\bm \te}} \exp\Big\{ \mathop{\sum\sum}_{j,k>0} \theta_{jk}S_j(u;X) S_k(v;Y)\Big\},\eeq
where $Z_{\bm \te}$ is the normalization factor that ensures $\cop_{\teb}$ is a proper density
\[ Z_{\bm \te}\,=\,\iint_{[0,1]^2} \exp\Big\{ \mathop{\sum\sum}_{j,k>0} \theta_{jk}S_j(u;X) S_k(v;Y)\Big\} \dd u \dd v.\]
We refer \eqref{exp:model} as the log-bilinear copula model.
\end{defn}
{\bf The Maximum Entropy Principle}.  Another justification for choosing the exponential model comes from the principle of maximum entropy (MaxEnt), pioneered by E. T. \cite{jaynes1957}. The maxent principle defines a unique probability distribution by maximizing the entropy 
$H(\cop)=-\int \cop_{X,Y} \log \cop_{X,Y}$ under the normalization constraint $\int \cop_{X,Y}=1$ and the following LP-co-moment conditions:
\beq \label{eq:cons}
\Ex_{\cop_{{\bm \te}}}[S_j(U;X) S_k(V;Y)]\,=\,\LP_{jk}.~\eeq

LP-co-means are orthogonal ``moments'' of copula, which can be estimated by 
\beq \label{eq:LPcom}
\tLP_{jk}\,=\,\Ex_{\widetilde \Cop}\big[ S_j(U;X) S_k(V;Y) \big]\,=\,\dfrac{1}{n} \sum_{i=1}^n S_j\big( \wtF_X(x_i)  ;X\big) \, S_k\big( \wtF_Y(y_i)  ;Y\big).\eeq
Applying calculus of variations, one can show that the maxent constrained optimization problem leads to the exponential \eqref{exp:model} form.  The usefulness of Jaynes' maximum entropy principle lies in providing a constructive mechanism to uniquely identify a probability distribution that is maximally non-committal (flattest possible) with regard to all unspecified information beyond the given constraints.

{\bf Estimation}. We fit a truncated exponential series estimator of copula density
\[\cop_{\teb}(u,v;X,Y)~=~\dfrac{1}{Z_{\bm \te}} \exp\Big\{ \sum_{j=1}^{m_1}\sum_{k=1}^{m_2} \theta_{jk}S_j(u;X) S_k(v;Y)\Big\}.\]
The task of finding the maximum likelihood estimates (MLE) of ${\bm \te}$ boils down to solving the following sets of equations for $j=1,\ldots,m_1$ and $k=1,\ldots,m_2$:
\beq \label{eq:mle}
\frac{\partial \log \cop_{\teb}}{\partial \te_{jk}} ~\equiv~ \frac{\partial \log Z_{\bm \te}}{\partial \te_{jk}}\,-\,\dfrac{1}{n} \sum_{i=1}^n S_j\big( \wtF_X(x_i)  ;X\big) \, S_k\big( \wtF_Y(y_i)  ;Y\big)\,=\,0.~~
\eeq
Note that the derivative of the log-partition function is equal to the expectation of the LP-co-mean functions:
\beq \label{eq:dlogp}
\dfrac{\partial \log Z_{\bm \te}}{\partial \te_{jk}}\,=\,\Ex_{\cop_{{\bm \te}}}[S_j(U;X) S_k(V;Y)].
\eeq
Replacing \eqref{eq:dlogp} and \eqref{eq:LPcom} into \eqref{eq:mle} implies that the MLE of MaxEnt model is same as the method of moments estimator satisfying the following moment conditions:
\[
\iint_{[0,1]^2} S_j(u;X) S_k(v;Y) \cop_{\teb}(u,v;X,Y)\dd u \dd v\,=\,\tLP_{jk}.\]
At this point, one can apply any convex optimization\footnote{since the second derivative of the log-partition function is a s positive semi-definite covariance matrix} routine (e.g., Newton's method, gradient descent, stochastic gradient descent, etc.) to solve for $\widehat{{\bm \te}}$.

\vskip.2em
{\bf Asymptotic}.  Let the sequence of $m_1$ and $m_2$ increase with sample size with an appropriate rate $\frac{(m_1 m_2 )^3}{n} \to 0$ as $\nti$. Then, under certain suitable regularity conditions, the exponential $\cop_{\widehat \teb}$ is a consistent estimate in the sense of Kullback-Leibler distance; see \cite{barron91} for more details.

\vskip.2em
{\bf Determining Informative Constraints}. 
Jayne's maximum entropy principle \textit{assumes} that a proper set of constraints (i.e., sufficient statistic functions) are given to the modeler, one that captures the phenomena under study. This assumption may be legitimate for studying thermodynamic experiments in statistical mechanics or for specifying prior distribution in Bayesian analysis, but certainly not for building empirical models. 

Which comes first: a parametric model or sufficient statistics? After all, the identification of significant components (sufficient statistics) is a prerequisite for constructing a legitimate probability model from the data \citep{D11a2}. Therefore the question of how to judiciously design and select the constraints from data, seems inescapable for nonparametrically learning maxent copula density function from data; also see Appendix \ref{app:VV}, which discusses the `two cultures' of maxent modeling. We address this issue as follows: (i) compute $\tLP_{jk}$ using the formula eq. \eqref{eq:LPcom}; (ii) sort them in  descending order based on their magnitude (absolute value); (iii) compute the penalized ordered sum of squares 
\[{\rm PenSum}(q)~=~\text{Sum of squares of top $q$ LP comeans}~-~\dfrac{\gamma_n}{n}q.~~\]
For AIC penalty choose $\gamma_n=2$, for BIC  choose $\gamma_n=\log n$, etc.  Further details can be found in \citet[Sec. 4.3]{D20copula}. (iv) Find the $q$ that maximizes the ${\rm PenSum}(q)$. Store the selected indices $(j,k)$ in the set $\mathcal{I}$. (v) Carry out maxent optimization routine based only on the selected LP-sufficient statistics-based constraints:
$$\Big\{ S_j(u;X)S_k(v;Y)\Big\},~~ (j,k) \in \mathcal{I}.~~~~~~~$$ 
This pruning strategy guards against overfitting. Finally, return the estimated \textit{reduced-order} (with effective dimension $|\mathcal{I}|$) maxent copula model. 

\begin{rem}[Nonparametric MaxEnt]
The proposed nonparametric maxent mechanism produces a copula density estimate, which is flexible (can adapt to the `shape of the data' without making risky a priori assumptions) and yet possesses a compact analytical form.
\vspace{-.4em}
\end{rem}

\subsection{Log-linear Model}
We provide a second parameterization of copula density. 
\begin{defn}
The log-linear orthogonal expansion of LP-copula is given by:
\beq \label{eq:lpspec} 
\cop_{\bm \mu}(u,v;X,Y)~=~\exp\big\{\mu_0 + \sum_{k>0} \mu_k \,\phi_k(u;X)\,\psi_k(v;Y)\big\},~~0<u,v<1,
\eeq
We call the parameters of this model ``log-linear LP-correlations'' that satisfy for $k>0$
\[
\mu_{k}=\iint_{[0,1]^2} \phi_k(u;X) \psi_k(v;Y) \log \cop_{\bm \mu}(u,v;X,Y) \dd u \dd v,~~
\]
\end{defn}
{\bf Connection}. Two fundamental representations, namely the log-bilinear  \eqref{exp:model} and loglinear \eqref{eq:lpspec} copula models, share some interesting connections\footnote{See \cite{D20copula} for a parallel result on the LP-orthogonal series copula model.}. To see that perform singular value decomposition (SVD) of the $\Theta$-matrix whose $(j,k)$th entry is $\te_{jk}$:
\[\Theta\,=\,U\Omega V^T\hspace{-.5em},\]
$u_{ij}$ and $v_{ij}$ are the elements of the singular vectors with singular values $\mu_1 \ge \mu_2 \ge \cdots \ge 0$. Then the spectral bases can be expressed as the linear combinations of the LP-polynomials:
\bea 
\phi_k(u;X)\,=\,\sum\nolimits_{j}u_{jk}S_{j}(u;X)~~\label{eq:sb1}\\
\psi_k(u;Y)\,=\,\sum\nolimits_{l}v_{lk}S_{l}(v;Y).~~\label{eq:sb2}
\eea
Hence, the LP-spectral functions (\ref{eq:sb1}-\ref{eq:sb2}) satisfy the following orthonormality conditions:
\beas 
\int\phi_k(u;X) \dd u\,=\,\int\psi_k(v;Y) \dd v\,=\,0~~~~~~\\
\hskip1em\int \phi_j(u;X) \psi_k(u;Y) \dd u\,=\,\delta_{jk}, \,\text{for}~ j\neq k.~~
\vspace{-.6em}
\eeas

\subsection{A Few Examples}
We demonstrate the flexibility of the LP-copula models using real data examples.
\vskip1.15em

\begin{figure}[ ]
  \centering
\includegraphics[width=.45\linewidth,keepaspectratio,trim=4cm 1.5cm 1cm 1cm]{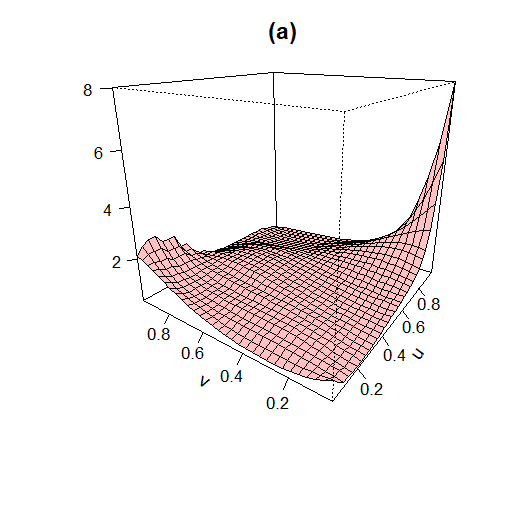}
\hspace{-.5em}
\includegraphics[width=.445\linewidth,keepaspectratio,trim=1cm 1.5cm 4cm 1cm]{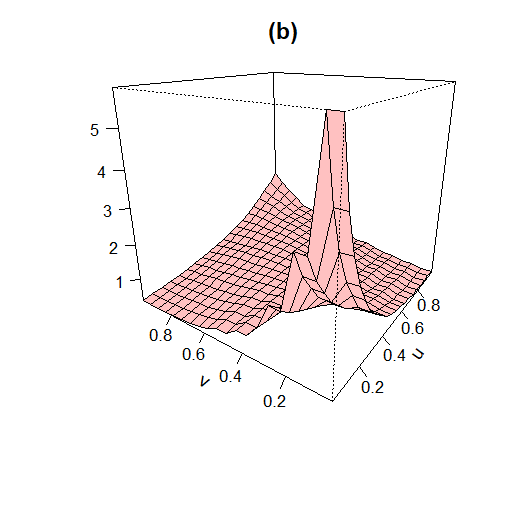}\\[2em]
\includegraphics[width=.45\linewidth,keepaspectratio,trim=4cm 1.5cm 1cm 1cm]{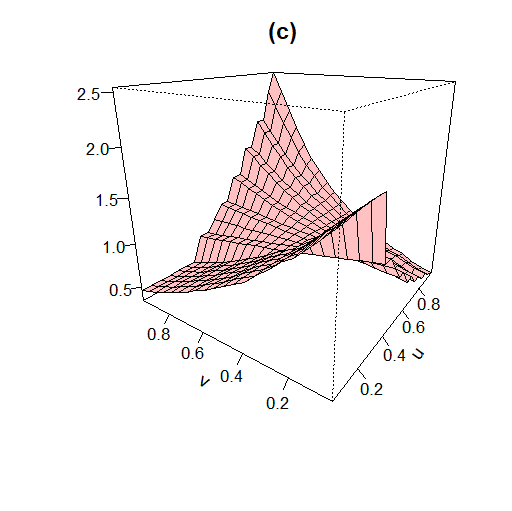}
\hspace{-.5em}
\includegraphics[width=.45\linewidth,keepaspectratio,trim=1cm 1.5cm 4cm 1cm]{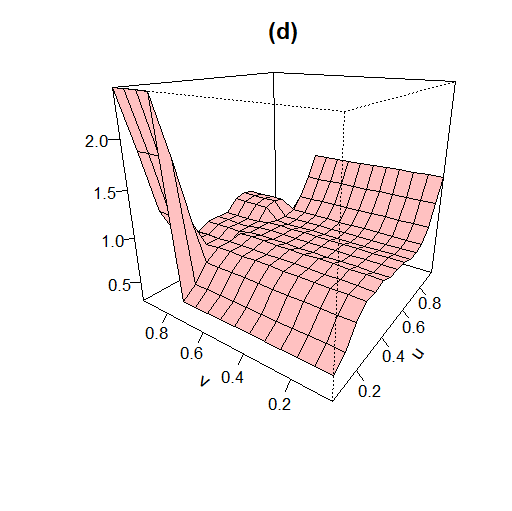}
\caption{MaxEnt LP-copula density estimates for (a) Kidney fitness (age vs tot: both continuous marginals), (b) PLOS data (title length vs number of authors: both discrete marginals), (c) horseshoe crab data (number of satellites vs width: mixed discrete-continuous marginals), and (d) Challenger space shuttle data (temperature vs number of O-ring failures: mixed continuous-discrete marginals).}
\label{fig:copFig1}
\vspace{-.25em}
\end{figure}

\begin{example} \label{ex:kidney} \textit{Kidney fitness data} \citep[Sec 1.1]{efron2016computer}.
It contains measurements on $n=157$ healthy volunteers (potential donors). For each volunteer, we have their age (in years) and a composite measure ``tot'' of overall function of kidney function.  To understand the relationship between age and tot, we estimate the copula: 
\[
\small{\whcop_{X,Y}(u,v)\,=\,\exp\big\{ -.40S_1(u;X)S_1(v;Y)+.18S_2(u;X)S_2(v;Y)-0.12  \big\},}
\]
displayed in Fig. \ref{fig:copFig1}(a). At the global scale, the shape of the copula density indicates a prominent negative ($\widehat \te_{11}=-0.40$) association between age and tot. Moreover, at the local scale, significant heterogeneity of the strength of dependence is clearly visible, as captured by the nonlinear asymmetric copula: the correlation between age and tot is quite high for older (say, $>70$) donors, compared to younger ones. This allows us to gain refined insights into \textit{how} kidney function declines with age.
\end{example}

\begin{example} \textit{PLOS data} (both discrete marginals). It contains information on $n=878$ journal articles published in PLOS Medicine between 2011 and 2015. For each article, two variables were extracted: length of the title and the number of authors. The dataset is available in the R-package \texttt{dobson}. The checkerboard-shaped estimated discrete copula
\[
\scalemath{0.95}{
\whcop_{X,Y}(u,v)\,=\,\exp\big\{ .42S_1(u;X)S_1(v;Y)+.10S_2(u;X)S_2(v;Y) -.07S_2(u;X)S_2(v;Y)-0.12  \big\}.
}
\] 
is shown in Fig. \ref{fig:copFig1}(b), which shows a strong positive nonlinear association. In particular, the sharp lower-tail, around the $(0,0)$, indicates that the smaller values of ($X,Y$) have a greater tendency to occur together than the larger ones.
\end{example}

\begin{example}\textit{Horseshoe Crabs Data} (mixed marginals).
The study consists of $n=173$ nesting horseshoe crabs \citep{agresti2014}. For each female crab in the study, we have its carapace width (cm) and number of male crabs residing nearby her nest. The goal of the study is to investigate whether carapace width affects number of male satellites for the female horseshoe crabs. If so, how--what is the shape of the copula dependence function?  The estimated copula, shown in Fig. \ref{fig:copFig1}(c), is given by:
\[
\small{\whcop_{X,Y}(u,v)\,=\,\exp\Big\{ 0.375S_1(u;X)S_1(v;Y)-0.077  \Big\}.}
\]
This indicates a significant positive linear correlation between the width and number of satellites of a female crab.
\end{example}

\begin{example} \textit{1986 Challenger Shuttle O-Ring data}. On January 28, 1986, just after seventy-three seconds into the flight, Challenger space shuttle broke apart, killing all seven crew members on board. The purpose of this study is to investigate whether the ambient temperature during the launch was related to the damage of shuttle's O-rings. For that we have $n=23$ previous shuttle missions data, consisting of  launch temperatures (degrees F), and number of damaged O-rings (out of 6). The estimated LP-maxent copula density
\[
\small{\whcop_{X,Y}(u,v)\,=\,\exp\Big\{ -0.37S_1(u;X)S_1(v;Y)-0.27S_3(u;X)S_1(v;Y)-0.14  \Big\}}
\]
is displayed in panel (d), and shows a strong negative association between the temperature at the launch and the number of damaged O-rings. Moreover, the sharp peak of the copula density around the edge $(0,1)$ further implies that cold temperatures can excessively increase the risk of failure of the o-rings. 
\vspace{-.5em}
\end{example}
\section{Applications to Statistical Modeling} \label{sec:Ut}
The scope of the general theory of the preceding section goes far beyond simply a tool for nonparametric copula approximation. In this section, we show how one can derive a large class of applied statistical methods in a unified manner by suitably reformulating them in terms of the LP-maxent copula model. In doing so, we also provide statistical interpretations of LP-maxent parameters under different data modeling tasks.


\subsection{Goodman's Association Model For Categorical Data}

Categorical data analysis will be viewed through the lens of LP-copula modeling. Let $X$ and $Y$ denote two discrete categorical variables; $X$ with $I$ categories and $Y$ with $J$ categories. The data are summarized in an $I\times J$ contingency table $\F$: $f_{kl}$ is the observed cell count in row $k$ and column $l$ of $\F$ and $n=f_{++}$ is the total frequency. The row and column totals are denoted as $f_{k+}$ and $f_{+l}$.  The observed joint $\Pr(X=k,Y=l)$ is denoted by $\tp_{kl}=f_{kl}/f_{++}$; the respective row and column marginals are given by $\tp_{k+}=f_{k+}/f_{++}$ and $\tp_{+l}=f_{+l}/f_{++}$. 

\vskip.5em
{\bf LP Log-linear Model}. We specialize our general copula model \eqref{eq:lpspec} for two-way contingency tables. The discrete LP-copula for the $I\times J$ table is given by
\beq \label{eq:copll}
\cop(F_X(k), F_Y(l)) \,=\, \exp\Big( \mu_0 + \sum_{j=1}^m \mu_{j} \phi_{jk}\psi_{jl}\Big),
\eeq 
where we abbreviate the row and columns scores $\phi_j(F_X(k))=\phi_{jk}$ and $\psi_j(F_Y(l))=\psi_{jl}$ for $k=1,2,\ldots,I$ and $l=1,2,\ldots,J$. The number of components $m \le M=\min(I-1,J-1)$; we call the log-linear model \eqref{eq:copll} `saturated' (or `dense') when we have $m=M$ components. The non-increasing sequence of model parameters $\mu_j$'s are called ``intrinsic association parameters'' that satisfy\footnote{Compare our equation \eqref{eq:parll} with  equation (34) of \citet{good96}.}
\beq \label{eq:parll}
\mu_j\,=\,\sum_{k=1}^I\sum_{l=1}^J \Big( \log p_{kl}\Big) p_{k+} p_{+l} \phi_{jk} \psi_{jl},~~\text{for}~j=1,\ldots,m.
\eeq 
Note that the discrete LP-row and column scores, by design, satisfy (for $j \neq j'$):
\bea
\sum_{k=1}^I\phi_{jk} p_{k+} \,=~\sum_{l=1}^J \psi_{jl}p_{+l} ~=~0 \label{eq:rs1}\\
\sum_{k=1}^I\phi^2_{jk} p_{k+} \,=~\sum_{l=1}^J \psi^2_{jl}p_{+l}~ =~1
\eea
\vspace{-2.2em}
\bea \label{eq:rs3}
\sum_{k=1}^I \phi_{jk} \phi_{j'k} p_{k+} \,=~\sum_{l=1}^J \phi_{jl} \psi_{j'l} p_{+l}\,=\,0.
\eea
\vskip.5em
{\bf Interpretation}. It is clear from \eqref{eq:parll} that the parameters $\mu_j$'s are fundamentally different from the standard Pearsonian-type correlation ${\rm Cor}(\phi_{j}, \psi_{j})=\Ex[\phi_{j}\psi_{j}]$, due to \eqref{eq:rs1}--\eqref{eq:rs3}:
\beq \label{eq:corP}
\rho_j\,=\,\sum_{k=1}^I\sum_{l=1}^J p_{kl}\phi_{jk} \psi_{jl},~~\text{for}\,j=1,\ldots,m.\eeq
The coefficients of the LP-MaxEnt-copula expansion for contingency tables carry a special interpretation in terms of log-odds-ratio. To see this we start by examining the $2\times 2$ case.
\vskip.5em
{\bf The $2\times 2$ Contingency Table}. Applying \eqref{eq:parll} for two-by-two tables we have
\beq \label{eq:mu1}
\mu_1\,=\,\sum_{k=0}^1\sum_{l=0}^1 \Big( \log p_{kl}\Big) p_{k+} p_{+l} \phi_{1k} \psi_{1l}.\eeq
Note that for dichotomous $X$ the LP-spectral basis $\phi_1(F_X(x))$ is equal to $T_1(x;F_X)$. Consequently, we have the following explicit formula for $\phi_1$ and $\psi_1$:
\bea 
\phi_1(F_X(x))&=&\dfrac{x-p_{2+}}{\sqrt{p_{1+} p_{2+}}}\label{eq:phibi1}\\
\psi_1(F_Y(y))&=&\dfrac{y-p_{+2}}{\sqrt{p_{+1} p_{+2}}} \label{eq:phibi2}
\eea
Substituting this into \eqref{eq:mu1} yields the following important result.
\begin{thm} \label{thm:2b2}
For 2-by-2 contingency tables, the estimate of the statistical parameter $\mu_1$ of the maxent LP-copula model 
\[ \cop_{\bm{\mu}}(u,v;X,Y)\,=\,e^{\mu_0 + \mu_1 \phi_1(u;X)\psi_1(v;Y)}\]
can be expressed as follows:
\beq \label{eq:mulor}
\widehat{\mu}_1\,=\,\log\Bigg[ \dfrac{\tp_{11}\tp_{22}}{ \tp_{12} \tp_{21}  } \Bigg] \Big(\tp_{1+}  \tp_{+1} \tp_{2+} \tp_{+2}\Big)^{1/2}\hspace{-.8em},
\eeq
where the part inside the square bracket is the sample log-odds-ratio.
\end{thm}
\begin{rem}[Significance of Theorem \ref{thm:2b2}]
We have \textit{derived} the log-odds-ratio statistic from first principles using a copula-theoretic framework.  To the best of our knowledge, no other study has discovered this connection; see also of \citet[eq. 16]{goodman1991} and \cite{gilula1988}. In fact, one can view Theorem \ref{thm:2b2} as a special case of the much more general result described next.
\end{rem}
\begin{thm}
For an $I\times J$ table, consider a two-by-two subtable with rows $k$ and $k'$ and columns $l$ and $l'$. Then the logarithm of odds-ratio $\vartheta_{kl,k'l'}$ is connected with the intrinsic association parameters $\mu_j$ in the following way:
\beq \label{eq:mulorG}
\log \vartheta_{kl,k'l'}\,=\,\sum_{j=1}^M \mu_j \big(\phi_{jk} - \phi_{jk'}\big)\big(\psi_{jl} - \psi_{jl'}\big).
\eeq
\end{thm}
To deduce \eqref{eq:mulor} from \eqref{eq:mulorG}, verify the following, utilizing the LP-basis formulae \eqref{eq:phibi1}-\eqref{eq:phibi2}
\beas 
\phi_{11}-\phi_{10}:=\phi_1(F_X(1)) - \phi_1(F_X(0))=\dfrac{p_{1+}+p_{2+}}{\sqrt{p_{1+}p_{2+}}} = \big(p_{1+}p_{2+}\big)^{-1/2}\\
~~\psi_{11}-\psi_{10}:=\psi_1(F_Y(1)) - \psi_1(F_Y(0))=\dfrac{p_{+1}+p_{+2}}{\sqrt{p_{+1}p_{+2}}} = \big(p_{+1}p_{+2}\big)^{-1/2} \hspace{-.4em}.
\eeas
\vskip.4em
{\bf Reproducing Goodman's Association Model}. Our discrete copula-based categorical data model \eqref{eq:copll} expresses the logarithm of ``dependence-ratios''\footnote{
\citet[pp. 410]{good96} calls it ``Pearson ratios.''}
\beq \label{eq:ratio}
\cop(F_X(k),F_Y(l))\,=\,\dfrac{p_{kl}}{p_{k+}p_{+l}}\eeq
as a linear combination of LP-orthonormal row and column scores satisfying \eqref{eq:rs1}-\eqref{eq:rs3}. The copula-dependence ratio \eqref{eq:ratio} measures the \textit{strength} of association between the $k$-th row category and $l$-the column category.  To make the connection even more explicit, rewrite \eqref{eq:copll} for two-way contingency tables as follows:
\beq \label{eq:Glogl}
\log p_{kl}\,=\,\mu_0 + \mu_k^{{\rm R}} + \mu_l^{{\rm C}} + \sum_{j=1}^m \mu_{j} \phi_{jk}\psi_{jl},
\eeq
where $\mu_k^{{\rm R}}$ denotes the logarithm of row marginal $\log p_{k+}$ and $\mu_l^{{\rm C}}$ denotes the logarithm of column marginal $\log p_{+l}$. \cite{goodman1991} called this  model \eqref{eq:Glogl} a ``weighted association model'' where weights are marginal row and column proportions. He used the term ``association model'' (to distinguish it from correlation \eqref{eq:corP} based model) as it studies the relationship between rows and columns using odds-ratio.
\begin{rem}
Log-linear models are a powerful statistical tool for categorical data analysis \citep{agresti2014}. Here we have provided a contemporary unified view of loglinear modeling for contingency tables from discrete LP-copula viewpoint. This newfound connection might open up new avenues of research.
\vspace{-.2em}
\end{rem}

\subsection{Logratio biplot: Graphical Exploratory Analysis}
We describe a graphical exploratory tool---logratio biplot, which allows a quick visual understanding of the relationship between the categorical variables $X$ and $Y$. In the following, we describe the process of constructing logratio biplot from the LP-copula model \eqref{eq:copll}.

\vskip.3em
{\bf Copula-based Algorithm}. Construct two scatter plots based on the top two dominant components of the LP-copula model: the first one is associated with the row categories, formed by the points $(\mu_1 \phi_{1k}, \mu_2 \phi_{2k})$ for $k=1,\ldots,I$; and the second one is associated with the column categories, formed by the points $(\mu_1 \psi_{1l}, \mu_2 \psi_{2l})$ for $l=1,\ldots,J$. Logratio biplot is a two-dimensional display obtained by overlaying these two scatter plots--the prefix `bi' refers to the fact that it shares a common set of axes for both the rows and columns categories.

\begin{figure}[ ]
  \centering
\includegraphics[width=.46\linewidth,keepaspectratio,trim=4cm 1.5cm .5cm 0cm]{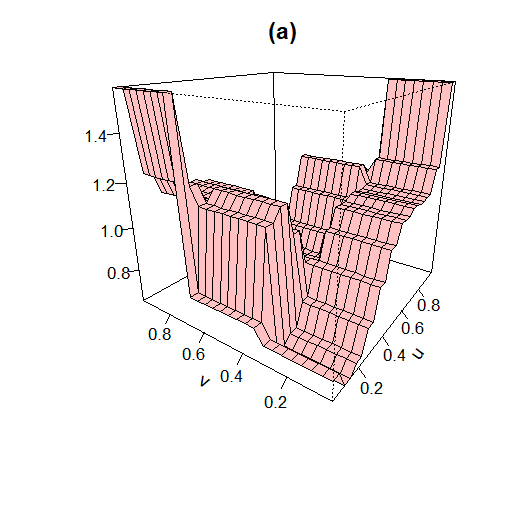}
\includegraphics[width=.44\linewidth,keepaspectratio,trim=1cm 0cm 3cm 5cm]{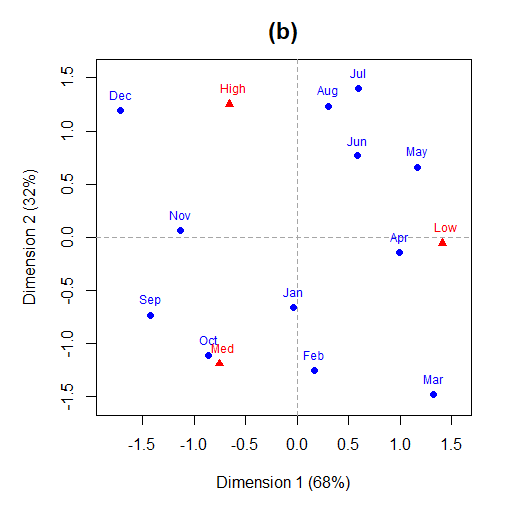}
\caption{1970 draft lottery: (a) piecewise-constant nonlinear LP-copula density estimate, and (b) two-dimensional logratio biplot that essentially captures all the useful information.}
\label{fig:draftUSA}
\end{figure}

\vskip.3em
{\bf Interpretation}. Here we offer an intuitive explanation of the logratio biplot from the copula perspective. We start by recalling the definition of conditional comparison density (CCD; see eq. \ref{eq:pbr}-\ref{eq:copC}), as the copula-slice. For fixed $X=k$, logratio-embedding coordinates $(\mu_j \phi_{jk})$ can be viewed as the LP-Fourier coefficients of the $d(F_Y(y); Y,Y|X=k)$, since
\[ 
d(F_Y(y); Y,Y|X=k)\,=\, \exp\Big\{ \mu_0 + \sum_{j=1}^m \big(\mu_{j} \phi_{jk}\big)\psi_{jl}\Big\},
\]
Similarly, the logratio coordinates $(\mu_j \psi_{jl})$ for fixed $Y=l$ can be interpreted as the LP-expansion coefficients of $d(F_X(x); X,X|Y=l)$. Hence, the logratio biplot can  alternatively be viewed as follows:  (i) estimate the discrete LP-copula density; (ii) Extract the copula slice $\whd(u; Y,Y|X=k)$ along with its LP-coefficients $(\mu_1 \phi_{1k}, \mu_2 \phi_{2k})$; (iii) similarly, get the estimated $\whd(v; X,X|Y=l)$---the copula slice at $Y=F_Y(l)$ along with its LP-coefficients $(\mu_1 \psi_{1l}, \mu_2 \psi_{2l})$; (iv) Hence, the logratio biplot (see Fig. \ref{fig:draftUSA}(b)) measures the association between the row and column categories $X=k$ and $Y=l$ by measuring the similarity between the `shapes' of $\whd(u; Y,Y|X=k)$ and $\whd(v; X,X|Y=l)$ through their LP-Fourier coefficients.

\vskip.25em
\begin{rem}[Historical Significance]
The following remarks are pertinent: (i) Log-ratio map traditionally taught and practiced using matrix-algebra \citep{greenacre2018compositional}. This is in sharp contrast with our approach, which has provided a \textit{statistical synthesis} of log-ratio biplot from a new copula-theoretic viewpoint. To the best of author's knowledge,  this is the first work that established such a connection. (ii) Logratio biplot has some important differences with the correspondence analysis pioneered by the French statistician Jean-Paul Benz{\'e}cri; for more details, see \cite{goodman1991} and \cite{Benzecri91}. However, in practice, these two methods often lead to very similar conclusions (e.g., contrast Fig. \ref{fig:draftUSA}(b) and  Fig. \ref{fig:CA}). 
\vskip.2em
\end{rem}

\begin{example} \textit{1970 Vietnam War-era US Draft Lottery} \label{ex:draft} \citep{fienberg1971rand}. All eligible American men aged 19 to 26 were drafted through a lottery system in 1970 to fill the needs of the country’s armed forces. In 1970, the US conducted a draft lottery to
determine the order (risk) of induction. The results of the draft are given to us in the form of a $12\times 3$ contingency table (see Table \ref{tab:draft} in the appendix): rows are months of the year from January to December, and columns denote three categories of risk of being drafted---high, medium, and low. The question is of interest whether
the lottery was fairly conducted; in other words, is there any association between the two categories of $12\times 3$ table? The discrete `staircase-shaped' LP-copula estimate is shown in the Fig. \ref{fig:draftUSA} (a), whose explicit form is given below:
\[
\whcop_{X,Y}(u,v;X,Y)~=~\exp\big\{0.26 \,\phi_1(u;X)\,\psi_1(v;Y)+0.18 \,\phi_2(u;X)\,\psi_2(v;Y)\,-\,0.052\big\}.
\]
We now overlay the scatter plots $(0.26\phi_{1k}, 0.18\phi_{2k})$ for $k=1,\ldots,12$ and $(0.26\psi_{1l}, 0.18\psi_{2l})$ for $l=1,2,3$ to construct the logratio biplot, as displayed in Fig. \ref{fig:draftUSA}. This easy-to-interpret two-dimensional graph captures the essential dependence pattern between the row (month: in blue dots) and the column (risk category: in red triangles) variables.
\end{example}

\subsection{Loglinear Modeling of Large
Sparse Contingency Tables} \label{sec:stable}
It has been known for a long time that
classical maximum likelihood-based log-linear models break down when applied to large sparse contingency tables with many zero cells; see \cite{fienberg20073c}. Here we discuss a new maxent copula-based smooth method for fitting a parsimonious log-linear model to sparse contingency tables.

\vskip.35em

\begin{example} \label{ex:zel} \textit{Zelterman data}. The dataset \cite[Table 1]{zelterman1987} is summarized as a $28\times 26$ cross-classified table that reports monthly salary and  number of years of experience since bachelor's degree of $n=129$ women employed as mathematicians or statisticians. The table is extremely sparse---86\% cells are empty! See Fig. \ref{fig:appzelterman2} of Appendix \ref{app:fig}.
\end{example}

\begin{figure}[ ]
  \centering
\includegraphics[width=.45\linewidth,keepaspectratio,trim=4cm 1cm 3cm 3cm]{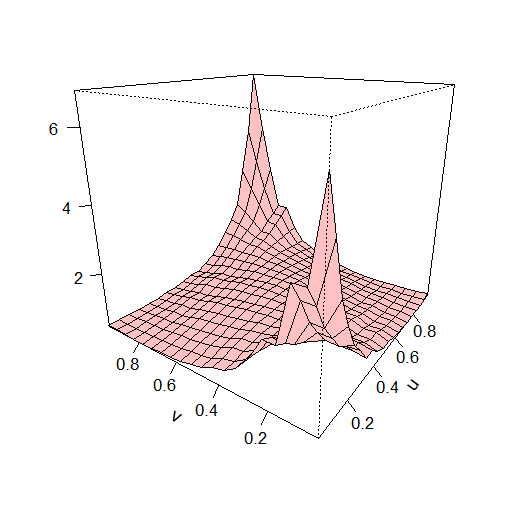}
\caption{Estimated LP-copula for the $28\times 26$ Zelterman contingency table data. A strong positive correlation is evident.}
\label{fig:zelterman}
\vspace{-.25em}
\end{figure}
\vskip.2em		
{\bf A Parsimonious Model}. The estimated smooth log-linear LP-copula model for the Zelterman data is given by:
\[
\scalemath{0.9}{
\whcop_{X,Y}(u,v)\,=\, \exp\Big\{0.52 \,\phi_1(u;X)\,\psi_1(v;Y)+0.37 \,\phi_2(u;X)\,\psi_2(v;Y)+0.18 \,\phi_3(u;X)\,\psi_3(v;Y)\,-\,0.27\Big\},}
\]
\vskip.1em
\vspace{-.15em}
displayed in Fig. \ref{fig:zelterman}. This shows a strong positive correlation between salary and number of years of experience. However, the most notable aspect is the effective model dimension, which can be viewed as the \textit{intrinsic} degrees of freedom (df). Our LP-maxent approach distills a compressed representation with reduced numbers of parameters that yields a \textit{smooth} estimates: only requiring $m=3$ components to capture the pattern in the data---a radical compression with negligible information loss! Contrast this with the dimension of the saturated loglinear model: $(28-1)\times (26-1)=675$ \textemdash a case of a severely overparameterized \textit{non-smooth} model with inflated degrees of freedom, which leads to an inaccurate goodness-of-fit test for checking independence between rows and columns. More on this in Sec. \ref{sec:G2}.

\vskip.3em		
{\bf Smoothing Ordered Contingency Tables}. The nonparametric maximum likelihood-based cell probability estimates $\tp_{X,Y}(k,l)=f_{kl}/n$ are very noisy and unreliable for sparse contingency tables. By sparse, we mean tables with a large number of cells relative to the number of observations.  

Using Sklar's representation theorem, one can simply estimate the joint probability $\hp_{X,Y}(k,l)$ by multiplying the empirical product pmf $\tp_X(k)\tp_Y(l)$ with the smoothed LP-copula. In other words, the copula can be viewed as a \textit{data-adaptive} bivariate discrete \textit{density-sharpening function} that corrects the independent product-density to estimate the cell probabilities. 
\beq~~~~~\texttt{dKernel}(k,l)\,=\,\whcop_{X,Y}\big(  \wtF_X(k), \wtF_Y(l)\big), ~~\text{for}~k=1,\ldots,I;~ l=1,\ldots,J.\eeq
where the discrete-kernel function satisfies
\[
\dfrac{1}{n^2} \sum_{k=1}^I\sum_{l=1}^J \texttt{dKernel}(k,l) f_{k+} f_{+l}\,=\,1.
\]
This approach can be generalized for \textit{any} bivariate discrete distribution; see next section.
\begin{rem}
For a comprehensive literature on traditional kernel-based nonparametric smoothing methods for sparse contingency tables,  readers are encouraged to consult \cite{simonoff1985,simonoff1995} and references therein.
\vspace{-.3em}
\end{rem}

\subsection{Modeling Bivariate Discrete Distributions}
The topic of nonparametric smoothing for multivariate discrete distributions has received far less attention than the continuous one.
Two significant contributions in this direction include: \cite{aitchison1976} and \cite{simonoff1983}. In what follows, we discuss a new  LP-copula-based procedure for  modeling correlated discrete random variables.

\vskip.35em
\begin{example} \label{ex:shun}
\textit{Shunter accident data} \citep{arbous1951accident}. As a motivating example, consider the following data: we are given the number of accidents incurred by $n=122$ shunters in two consecutive year periods, namely 1937-1942 and 1943-1947. To save space, we display the bivariate discrete data in a contingency table format; see Table \ref{tab:shunter}.
\end{example}

\textit{Algorithm}.  The main steps of our analysis are described below.

\begin{figure}[ ]
  \centering
\includegraphics[width=.454\linewidth,keepaspectratio,trim=1cm 1cm 1cm 1cm]{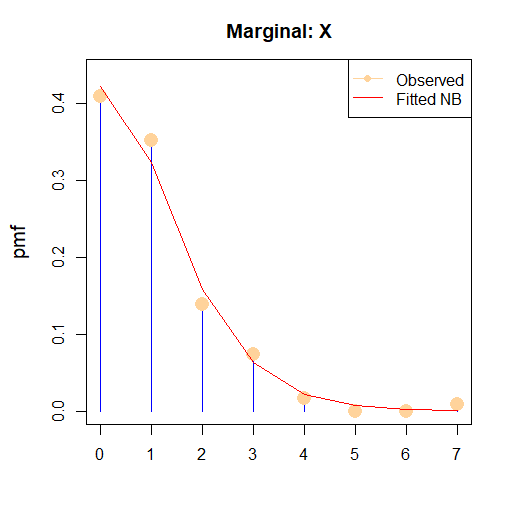}~~~~~~ \includegraphics[width=.454\linewidth,keepaspectratio,trim=1cm 1cm 1cm 1cm]{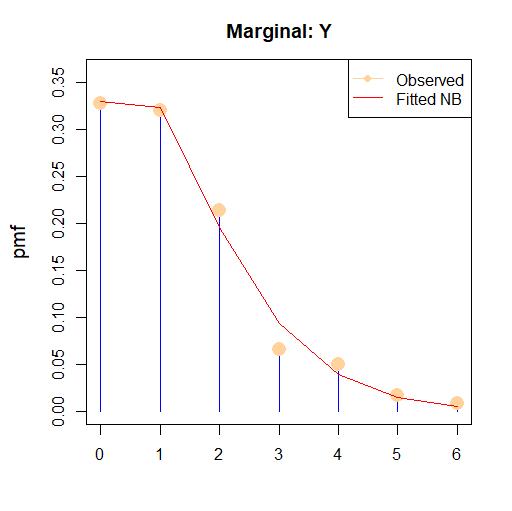}\\[3.55em]
\includegraphics[width=.6\linewidth,keepaspectratio,trim=3cm 1cm 3cm 2cm]{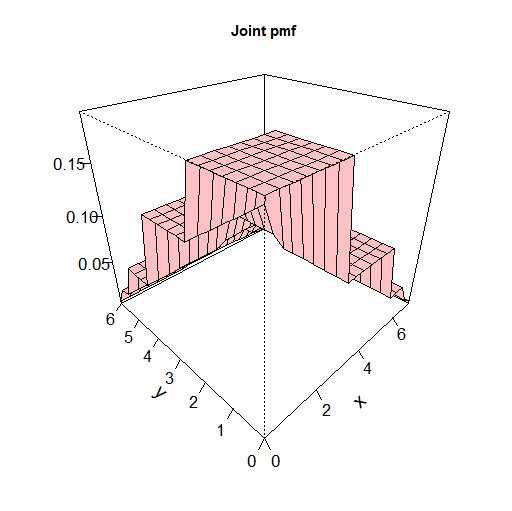}
\vskip1em
\caption{Accident Data. Top panel: marginal modeling---comparing the observed empirical pmf with the fitted negative binomial (NB) distribution. Bottom panel shows the estimated smooth joint pmf, which is obtained by ``sharpening'' the product of marginal densities (top row) using the LP-smooth copula (Eq. \ref{eq:shunCOP}).}\label{fig:ACC}
\end{figure}

~\texttt{Step 1}. \textit{Modeling marginal distributions}. We start by looking at the marginal distributions of $X$ and $Y$. As seen in Fig. \ref{fig:ACC}, negative binomial distributions provide excellent fit. To fix the notation, by $G_{\mu,\phi}= {\rm NB}(y;\mu,\phi)$, we mean the following probability distribution:
\[ {\rm NB}(y;\mu,\phi) = \binom{y + \phi - 1}{y} \,
\left( \frac{\mu}{\mu+\phi} \right)^{\!y} \, \left(
\frac{\phi}{\mu+\phi} \right)^{\!\phi} \!,~~y \in \mathbb{N},\]
where $\Ex(X)=\mu$ and $\Var(X)=\mu+\frac{\mu^2}{\phi}$. Using the method of MLE, we get: $X \sim G_1 ={\rm NB}(x;\hat \mu=0.97,\hat \phi=3.60)$, and $Y \sim G_2 ={\rm NB}(y; \hat \mu=4.30,\hat \phi=1.27)$.

\begin{table}[t]
\vskip.65em
\setlength{\tabcolsep}{11pt}
\renewcommand{\arraystretch}{.88}
\centering
\begin{tabular}{c|rrrrrrr}
  \hline
  & \multicolumn{7}{c}{1943-47} \\
    \cline{2-8}
1937-42& 0 &1 &2 &3 &4 &5 &6 \\
 \hline
0 &  21 &  18 &   8 &   2 &   1 &   0 &   0 \\ 
  1 &  13 &  14 &  10 &   1 &   4 &   1 &   0 \\ 
  2 &   4 &   5 &   4 &   2 &   1 &   0 &   1 \\ 
  3 &   2 &   1 &   3 &   2 &   0 &   1 &   0 \\ 
  4 &   0 &   0 &   1 &   1 &   0 &   0 &   0 \\ 
  7 &   0 &   1 &   0 &   0 &   0 &   0 &   0 \\ 
   \hline
\end{tabular}
\vskip.5em
\caption{Shunter accidents data displayed compactly as XY-contingency table.} \label{tab:shunter}
\end{table}
\vskip.3em
~\texttt{Step 2}. \textit{Generalized copula density}. 
The probability of bivariate distribution at $(x,y)$ can be written as follows (generalizing Sklar's Theorem): 
\beq \label{eq:gsk} 
\Pr(X=x,Y=y):=p_{X,Y}(x,y)=g_1(x) g_2(y) \,\cop_{X,Y}\big(G_1(x), G_2(y)\big),
\eeq
where generalized log-copula density admits the following decomposition:
\beq 
\log \big\{\hspace{-.15em}\cop_{X,Y}(G_1(x), G_2(y))\big\}\,=\, \sum_{j=1}^{m_1}\sum_{k=1}^{m_2} \theta_{jk}T_j(x;G_1) T_k(y;G_2) \,-\,\log Z_{\bm \te}\,.
\eeq
It is important to note that the set of LP-basis functions $\{T_j(x;G_1)\}$ and $\{T_k(y;G_2)\}$ are specially designed for the parametric marginals $G_1$ and $G_2$, obeying the following \textit{weighted} orthonormality conditions:  
\beas \sum\nolimits_x g_1(x) T_j(x;G_1)=0,&\text{and}&\sum\nolimits_x g_1(x) T_j(x;G_1)T_k(x;G_1)=\delta_{jk}; \\
\sum\nolimits_x g_2(x)T_j(x;G_2)=0,&\text{and}&\sum\nolimits_x g_2(x) T_j(x;G_2)T_k(x;G_2)=\delta_{jk}.
\eeas
We call them  \texttt{gLP}-basis, to distinguish them from the earlier empirical LP-polynomial systems $\{T_j(x;\wtF_X)\}$ and $\{T_k(y;\wtF_Y)\}$; see Appendix \ref{app:LPb}.

\begin{rem}[Generalized copula as density-sharpening function]
The generalized copula
\beq \label{sup:dsf2d1}
\cop_{X,Y}\big(G_1(x), G_2(y)\big) = \dfrac{p_{X,Y}(x,y)}{g_1(x) g_2(y)},\eeq
acts as a bivariate ``density sharpening function'' in \eqref{eq:gsk} that \textit{corrects} the possibly \textit{misspecified} $g_1(x) g_2(y)$. This is very much in the spirit of \cite{D21discrete,D22ADM}. It is also instructive to contrast our generalized copula \eqref{sup:dsf2d1} with the usual definition of copula (c.f. Sec. \ref{sec:back}):
\[
\cop_{X,Y}\big(F_X(x), F_Y(y)\big) = \dfrac{p_{X,Y}(x,y)}{p_X(x) \,p_Y(y)}, \]
which requires correct specification of the marginals $p_X(x)$ and $p_Y(y)$.
\end{rem}
~\texttt{Step 3}. \textit{Exploratory goodness-of-fit.} The estimated log-bilinear LP-copula is
\beq \label{eq:shunCOP}
\whcop_{X,Y}(u,v)\,=\,\exp\Big\{ 0.287S_1(u;G_1)S_1(v;G_2)-0.043  \Big\}.
\eeq
There are three important conclusions that can be drawn from this non-uniform copula density estimate: (i) Goodness of fit diagnostic: the independence model (product of parametric marginals) $g_\perp(x,y)=g_1(x) g_2(y)$ is not adequate for the data. (ii) Nature of discrepancy: the presence of significant $\hat\te_{11}=0.287$ in the model \eqref{eq:shunCOP} implies that the tentative independence model should be updated by incorporating the strong (positive) `linear' correlation between $X$ and $Y$. (iii) Nonparametric repair: how to \textit{update} the initial $g_\perp(x,y)$ to construct a ``better'' model? Eq. \eqref{eq:gsk} gives the general updating rule, which simply says: copula provides the necessary bivariate-correction function to reduce the `gap' between the starting misspecified model $g_\perp(x,y)$ and the true unknown distribution $p_{X,Y}(x,y)$. (iv) In contrast to unsmoothed empirical multilinear copulas \citep{genest2013b}, our method produces smoothed and compactly parametrizable $\whcop(u,v)$ for discrete data.

\vskip.3em
~\texttt{Step 4}. \textit{LP-smoothed probability estimation.}  The bottom panel Fig. \ref{fig:ACC} shows the final smooth probability estimate $\hp_{X,Y}(x,y)$, computed by substituting \eqref{eq:shunCOP} into \eqref{eq:gsk}. Also compare Tables \ref{tab:accEMP} and \ref{tab:accLP}
of Appendix \ref{app:fig}.

\begin{rem}
Our procedure fits a `hybrid' model: a nonparametrically corrected (through copula) multivariate parametric density estimate.\footnote{A similar philosophy was proposed in \cite{Deep17LPMode} for univariate continuous distribution case.}  
One can use any parametric distribution instead of a negative binomial. The algorithm remains fully automatic, irrespective of the choice of parametric marginals $G_1$ and $G_2$, which makes it a universal procedure.
\end{rem}

\subsection{Mutual Information}
Mutual information (MI) is a fundamental quantity in Statistics and Machine Learning, with wide-ranging applications from neuroscience to physics to biology. For continuous random variables $(X,Y)$, mutual information is defined as
\beq \label{eq:MI0}
\MI(X,Y)\,=\,\iint f_{X,Y}(x,y) \log \dfrac{ f_{X,Y}(x,y) }{f_X(x) f_Y(y)} \dd x \dd y. 
\eeq
Among non-parametric MI estimators, $k$-nearest-neighbor and kernel-density-based methods \citep{moon1995MI,kraskov2004estimating, zeng2018jackknife} are undoubtedly the most popular ones. Here we are concerned with a slightly general problem of developing a flexible MI estimation algorithm that is: (D1) applicable for mixed\footnote{Reliably estimating MI for mixed case is notoriously challenging task \citep{gao2017estimating}.} $(X,Y)$; (D2) robust in the presence of noise; and,
(D3) invariant under monotone transformations\footnote{This is essential to make the analysis less sensitive to various types of data preprocessing, which is done routinely in applications like bioinformatics, astronomy, and neuroscience.}. To achieve this goal, we start by rewriting MI \eqref{eq:MI0} using copula:
\beq  \label{eq:cMI}
\MI(X,Y)\,=\,\int_{[0,1]^2} \cop_{X,Y}(u,v) \log \cop_{X,Y}(u,v) \dd u \dd v.
\eeq 

The next theorem presents an elegant closed-form expression for MI in terms of LP-copula parameters, which allows a fast and efficient estimation algorithm.
\vskip1em
\begin{thm} Let $(X,Y)$ be a mixed-pair of random variables. Under the LP log-bilinear copula model \eqref{exp:model}, the mutual information between $X$ and $Y$ has the following representation in terms of LP-co-mean parameters $\LP_{jk}$ and maximum entropy coefficients $\te_{jk}$
\beq
\MI_{\teb}(X,Y)\,=\,\mathop{\sum\sum}_{j,k>0} \te_{jk}\LP_{jk}\, -\, \log Z_{\teb}\,.
\eeq
\end{thm}
\vspace{-.35em}
{\bf Proof}.\, Express mutual information as:
\[
\MI_{\teb}(X,Y)\,=\,\Ex_{X,Y}\big[ \log \cop_{\teb}\big]\,=\,\sum_j\sum_k \te_{jk} \Ex_{X,Y}\big[ S_j(U;X) S_k(V;Y)\big]\,-\,\log Z_{\teb}\,.\]
The first equality follows from \eqref{eq:cMI} and the second one from \eqref{exp:model}. Complete the proof by replacing by $\LP_{jk}$ by $\Ex[ S_j(U;X) S_k(V;Y)]$ by virtue of \eqref{eq:cons}. As a practical consequence, we have the following efficient and direct MI-estimator, satisfying D1-D3:
\beq \label{eq:MIest}
\widehat{\rm MI}_{\teb}(X,Y)\,=\,\mathop{\sum\sum}_{j,k>0} \hte_{jk}{\tLP}_{jk}\, -\, \log Z_{\widehat \teb}\,.
\vspace{-.5em}
\eeq
\vskip.2em
{\bf Bootstrap inference.} Bootstrap provides a convenient way to estimate the standard error of the estimate \eqref{eq:MIest}. Perform bootstrap sampling, i.e., sample $n$ pairs of $(x_i,y_i)$ with replacement and compute $\widehat{\rm MI}$. Repeat the process, say, $B=500$ times to get the sampling distribution of the statistic. Finally, return the standard error of the bootstrap sampling distribution along with 95\% percentile-confidence interval.

\textit{Continuous $(X,Y)$ example}. Consider the kidney fitness data, discussed in Example \ref{ex:kidney}. The LP-copula-based (using $m=4$) method yields: $\hMI=0.230 \,\,(\pm 0.021)$. To understand how precise is the estimate, we have reported the bootstrap standard error in parentheses.

\begin{rem}
MI \eqref{eq:cMI} measures the departure of copula density from uniformity. This is because, MI can be viewed as the Kullback-Leibler (KL) divergence between copula and the uniform density: ${\rm MI}(X, Y)={\rm KL}(\cop; U_{[0,1]^2})$. A few immediate consequences: (i) MI is always nonnegative, i.e., ${\rm MI}(X, Y) \ge 0$,
and equality holds if and only if variables are independent. Moreover, the stronger the dependence between two variables, the larger the MI. (ii) MI is also invariant under different marginalizations. Two additional applications of MI (for categorical data and feature selection problems) are presented below. 
\vspace{-.4em}
\end{rem}

\subsubsection{Application 1:~\,(X,Y) Discrete: Smooth-G$^2$ Statistic} \label{sec:G2}
Given $n$ independent samples from an $I\times J$ contingency table, the $G^2$-test of goodness-of-fit, also known as the log-likelihood ratio test\footnote{In 1935, Samuel Wilks introduced log-likelihood ratio test as an alternative to Pearson’s chi-square test. In our notation, Pearson proposed $\int \cop^2_{X,Y}(u,v)$ and Wilks proposed $2\times \int \cop_{X,Y}(u,v)\log \cop_{X,Y}$---both are conceptually  equivalent: measuring how much the copula density deviates from the uniformity.}, is defined as
\beq \label{eq:G2}
G^2(X,Y)\,=\,2n\sum_{k=1}^I\sum_{l=1}^J \,\tp_{kl} \log \dfrac{\tp_{kl}}{\tp_{k+}\tp_{+l}},
\eeq
which under the null hypothesis of independence has asymptotic $\chi^2_{(I-1)(J-1)}$ distribution. From \eqref{eq:G2} one can immediately conclude the following.
\begin{thm}
The $G^2$ log-likelihood ratio statistic can be viewed as the raw nonparametric MI-estimate
\beq \label{eq:MIG}
G^2(X,Y)/2n\,=\,\tMI(X,Y),~~~~
\eeq 
where $\tMI(X,Y)$ is obtained by replacing the unknown distributions in \eqref{eq:MI0} with their empirical estimates.
\end{thm}
\begin{example} \textit{Hellman's Infant Data} \citep{yates1934} We verify the identity \eqref{eq:MIG} for the following $2\times 2$ table \ref{tab:yates}, which shows cross-tabulation of $n=42$ infants based on whether the infant was breast-fed or bottle-fed. 
\begin{table}[h]
\vskip1em
\begin{center}
\setlength{\tabcolsep}{1.15em}
\begin{tabular}{c|cc}
& Normal teeth & Malocclusion\\
\hline\\[-.4em]
Breast-fed & 4 & 16 \\[.2em]
Bottle-fed & 1& 21\\[.2em]
\hline
\end{tabular}
\end{center}
\caption{The data table on malocclusion of the teeth in infants were obtained by M. Hellman and reported in the classic paper by Frank Yates (1934, p.230).} \label{tab:yates}
\vspace{-.7em}
\end{table}

The estimated LP-copula density (shown in Fig. \ref{fig:hellmapp} of Appendix \ref{app:fig}) is given  by
\beq 
\whcop(u,v;X,Y)=\exp\big\{ 0.234S_1(u;X)S_1(v;Y)\,-\,0.03 \big\}.
\eeq
The empirical MI estimate is given by $2n \times  \tMI(X,Y)=2.50$, with the pvalue $0.12$ computed using the asymptotic null distribution $\chi^2_1$. This exactly matches with the $G^2$-statistic value; one may use the R-function \texttt{GTest}. 
\end{example}

The problem arises when we try to apply $G^2$-test for large sparse tables, and it is not hard to see why: the adequacy of asymptotic $\chi^2_{(I-1)(J-1)}$ distribution depends both on the sample size $n$ and the number of cells $p = IJ$. \cite{koehler1986goodness} showed that the approximation completely breaks down when $n/p<5$, leading to erroneous statistical inference due to significant loss of power; see Appendix \ref{app:simulation}. The following example demonstrates this.

\vskip.3em

\begin{example}
\textit{Zelterman Data Continued}. Log-likelihood ratio $G^2$-test produces pvalue $1$, firmly concluding the independence between salary and years of experience. This directly contradicts our analysis of Sec. \ref{sec:stable}, where we found a clear positive dependence between these two variables. Why $G^2$-test was unable to detect that effect? Because it is based on chi-square approximation with degrees of freedom $(28-1)\times (26-1)=675$. This inflated degrees of freedom completely ruined the power of the test. To address this problem, we recommend the following \textit{smoothed} version:
\beq \label{eq:sgMI}
{\rm Smooth-}G^2(X,Y)/2n=\hMI(X,Y),~~~~~
\eeq
where $\hMI(X,Y)$ is computed based on the LP-bilinear copula model:
\[
\whcop(u,v;Y,X)=\exp\Big\{0.50 \,S_1(u;Y)\,S_1(v;X)+0.18 \,S_2(u;Y)S_2(v;X)\,-\,0.167\Big\},
\]
${\rm Smooth-}G^2$ analysis (with df$=2$) generates pvalue $2.48\times 10^{-11}$, thereby successfully detecting the association.  The crucial aspect of our approach lies in its ability to provide a reduced dimensional parametrization of copula density. For the Zelterman data, we need just two components (i.e., the effective degrees of freedom is $2$) to capture the pattern.  
\end{example}

\begin{rem}[Discrete variables with many
categories]
Discrete distributions over large domains routinely arise in large-scale biomedical data such as diagnosis codes, drug compounds and genotypes \citep{seok2015MI}. The method proposed here can be used to jointly model such random variables. 
\end{rem}

\subsubsection{Application 2:~\,(X,Y) Mixed: Feature Importance Score} 
We consider the two-sample feature selection problem where $Y$ is a binary response variable, and $X$ is a predictor variable that can be either discrete or continuous.

\begin{figure}[ ]
  \centering
  \includegraphics[width=.57\linewidth,keepaspectratio,trim=3cm 3cm 3cm 4cm]{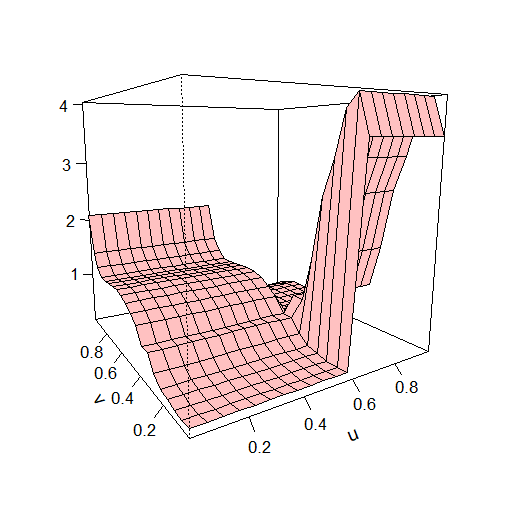}\\[1.36in]
\includegraphics[width=.424\linewidth,keepaspectratio,trim=3cm 2cm 1cm 4cm]{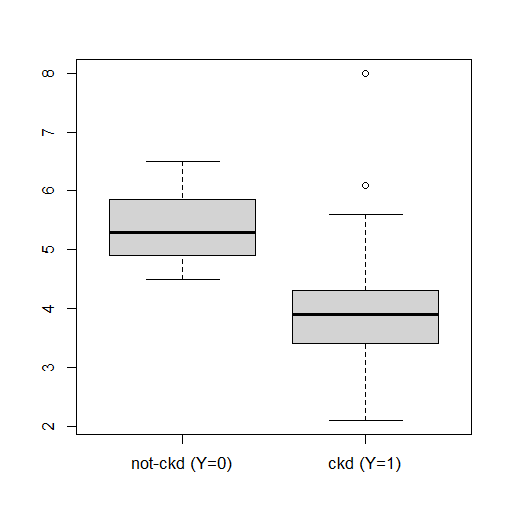}~~~~~~
\includegraphics[width=.424\linewidth,keepaspectratio,trim=1cm 2cm 3cm 4cm]{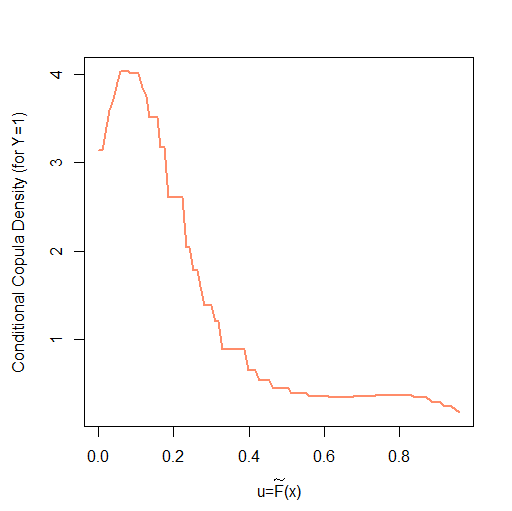}
\vskip2.5em
\caption{CKD data: The top panel shows the estimated copula density. The bottom panel shows the two-sample boxplots and the conditional comparison density (CCD) $d(v;X,X|Y=1)$. The piecewise constant shape of estimated CCD is not an aberration of our nonparametric approximation method; it reflects the inherent discreteness of the feature $X$ (red blood cell count).}
\label{fig:ckd1}
\vspace{-.25em}
\end{figure}

\begin{example} \textit{Chronic Kidney Disease data}. The goal of this study is to investigate whether kidney function is related to red blood cell count (RBC). We have a sample of $n=203$ participants, among whom $79$ have chronic kidney disease (ckd) and another $124$ are non-ckd. $Y$ denotes the kidney disease status and $X$ denotes the measurements on RBC (unit in million cells per cubic millimeter of blood). The estimated LP-copula density 
\beq 
\label{eq:ckdcop} 
\scalemath{0.82}{
\whcop(u,v)=\exp\Big\{-0.76 \,S_1(u;Y)\,S_1(v;X)+0.18 \,S_1(u;Y)S_2(v;X)-0.19 \,S_1(u;Y)S_4(v;X)\,-\,0.33\Big\},}
\eeq
is shown in Fig. \ref{fig:ckd1}.  We propose mutual information-based feature importance measure based on the formula \eqref{eq:MIest}; this yields $\hMI(Y,X)=0.36$ with pvalue almost zero, strongly indicating that RBC is an importance risk-factor related to kidney dysfunction.

What additional insights can we glean from this copula? To answer that question, let's focus our attention on the copula-slice for $u \in [0.61,1]$. This segment of the copula is essentially the conditional comparison density $d(v;X,X|Y=1)$ (see Sec. \ref{sec:back}), which can be easily derived by substituting $S_1(\wtF_Y(1);Y)=\sqrt{\frac{1-\tilde{\mu}}{\tilde{\mu}}}=1.25$ into \eqref{eq:ckdcop}, where $\tilde{\mu}=79/203=0.389$:
\[
\widehat d(v;X,X|Y=1)\,=\,\exp\big\{ -0.95S_1(v;X) + 0.23 S_2(v;X)-0.24S_4(v;X) -0.33\big\}.
\]
\vspace{-2em}
\end{example}

A few remarks on the interpretation of the above formula: 

~$\bullet$ \textit{Distributional effect-size}: Bearing in mind Eqs. (\ref{eq:du}, \ref{eq:copC}), note that $d(v;X,X|Y=1)$ compares two densities: $f_{X|Y=1}(x)$ with $f_X(x)$, thereby capturing the distributional difference. This has advantages over traditional two-sample feature importance statistic (e.g., Student's t or Wilcoxon statistic) that can only measure differences in location or mean.

\vskip.3em
~$\bullet$ \textit{Explainability}: The estimated $\whd(v;X,X|Y=1)$ involves three significant LP-components of $X$; the presence of the 1st order `linear' $S_1(v;X)$ indicates location-difference; the 2nd order `quadratic'  $S_2(v;X)$ indicates scale-difference; and 4th order `quartic' $S_4(v;X)$ indicates the presence of tail-difference in the two RBC-distributions. In addition, the negative sign of the linear effect $\hte_{11}=-0.95$ implies reduced mean level of RBC in the ckd-population. Medically, this makes complete sense, since a dysfunctional kidney cannot produce enough Erythropoietin (EPO) hormone, which causes the RBC to drop.


\subsection{Nonparametric Copula-Logistic Regression}
We describe a new copula-based nonparametric logistic regression model.  The key result is given by the following theorem, which provides a first-principle derivation of a robust nonlinear generalization of the classical linear logistic regression model.
\begin{thm} \label{tmh:logis} Let $\mu=\Pr(Y=1)$ and $\mu(x)=\Pr(Y=1|X=x)$. Then we have the following expression for the logit (log-odds)  probability model:
\beq \label{coplogthm}
\logit \left\{\mu(x)\right\} \,=\, \log \left( \dfrac{\mu(x)}{1-\mu(x)}\right)\,=\, \al_0 \,+\sum_j \al_j T_j(x;F_X),~~~~~\eeq
where $\al_0=\logit(\mu)$ and $\al_j=\frac{\te_{j1}}{\sqrt{\mu(1-\mu)}}$.
\end{thm}
{\bf Proof}. The proof consists of four main steps.  

Step 1. To begin with, notice that for $Y$ binary and $X$ continuous, the general log-bilinear LP-copula density function \eqref{exp:model} reduces to the following form: 
\beq  \label{eq:copmixed}
\cop_{\teb}(u,v;X,Y)~=~\dfrac{1}{Z_\te} \exp\Big\{ \sum_j \te_{j1} S_j(u;X) S_1(v;Y)\Big\},
\eeq
since we can construct at most $2-1=1$ LP-basis function for binary $Y$.

Step 2. Apply copula-based Bayes Theorem (Eq. \ref{eq:copG}) and express the conditional comparison densities as follows:
\beq \label{eq:bd1}
d_1(x)\,\equiv\, d(F_X(x);X, X|Y=1)\,=\,\dfrac{\mu(x)}{\mu},~~~~~
\eeq
and also,
\beq \label{eq:bd0}
d_0(x)\,\equiv\, d(F_X(x);X, X|Y=0)\,=\,\dfrac{1-\mu(x)}{1-\mu}.
\eeq
Taking logarithm of the ratio of \eqref{eq:bd1} and \eqref{eq:bd0}, we get the following important identity:
\beq \label{eq:mainLOG}
\log \left( \dfrac{\mu(x)}{1-\mu(x)}\right)\,=\,\log \left( \dfrac{\mu}{1-\mu}\right)\,+\,\log d_1(x)\,-\,\log d_0(x).
\eeq

Step 3. From \eqref{eq:copmixed}, one can deduce the following orthonormal expansion  of maxent-conditional copula slices $\log d_1$ and $\log d_0$:
\bea 
\log d_1(x)&=&  \sum_j\big( \te_{j1} T_1(1;F_Y)  \big)\, T_j(x;F_X) \,-\,\log Z_\te \label{eq:logd1}\\
\log d_0(x)&=& \sum_j\big( \te_{j1} T_1(0;F_Y)  \big) \, T_j(x;F_X)  \,-\,\log Z_\te \label{eq:logd0}
\eea 

Step 4. Substituting \eqref{eq:logd1} and \eqref{eq:logd0} into \eqref{eq:mainLOG} we get:
\[
\log \left( \dfrac{\mu(x)}{1-\mu(x)}\right)\,=\,\log \left( \dfrac{\mu}{1-\mu}\right)\,+\,\sum_j \Big\{ \frac{\te_{j1}}{\sqrt{\mu(1-\mu)}}\Big\} T_j(x;F_X),
\]
since for binary $Y$ we have (see Appendix \ref{app:LPb}):
\[
T_1(1;F_Y)-T_1(0;F_Y)\,=\,\dfrac{1-\mu}{\sqrt{\mu(1-\mu)}}\,+\,\dfrac{\mu}{\sqrt{\mu(1-\mu)}}\,=\,\dfrac{1}{\sqrt{\mu(1-\mu)}}. 
\]
Substitute $\al_0=\logit(\mu)$ and $\al_j=\frac{\te_{j1}}{\sqrt{\mu(1-\mu)}}$ to complete the proof. \qed

\begin{figure}[ ]
  \centering
  \includegraphics[width=.44\linewidth,keepaspectratio,trim=2cm 2cm 1cm 1cm]{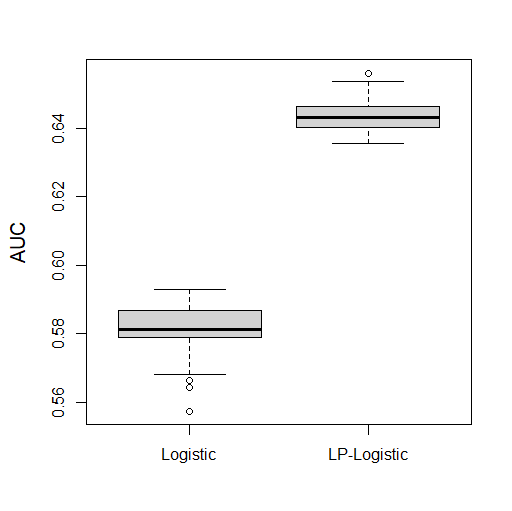}~~~~
  \includegraphics[width=.44\linewidth,keepaspectratio,trim=1cm 2cm 2cm 1cm]{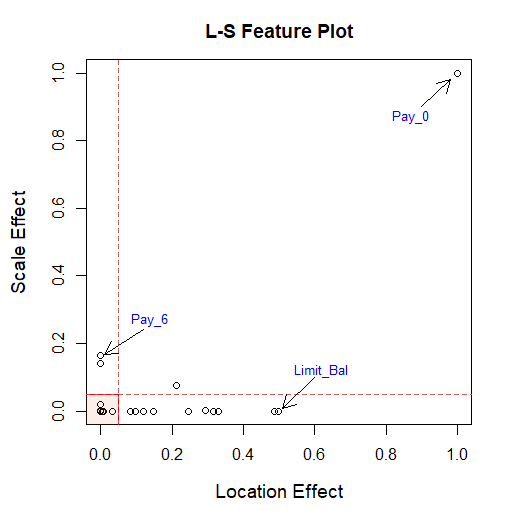}\\[5.5em]
 \includegraphics[width=.32\linewidth,keepaspectratio,trim=.5cm 1cm .5cm .5cm]{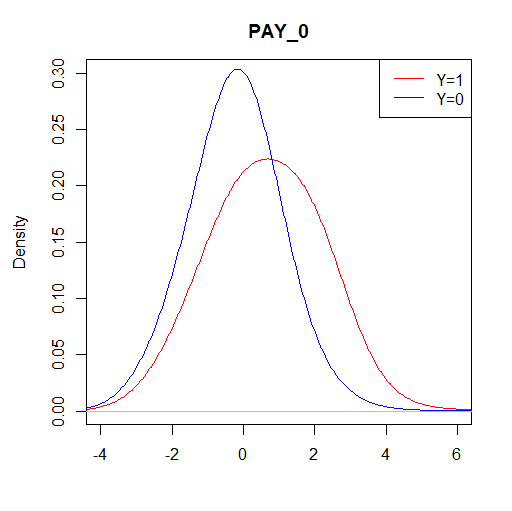} ~
   \includegraphics[width=.32\linewidth,keepaspectratio,trim=.5cm 1cm .5cm .5cm]{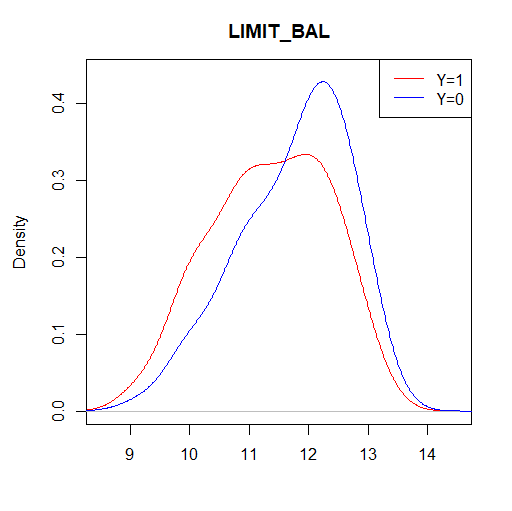}~
    \includegraphics[width=.32\linewidth,keepaspectratio,trim=.5cm 1cm .5cm .5cm]{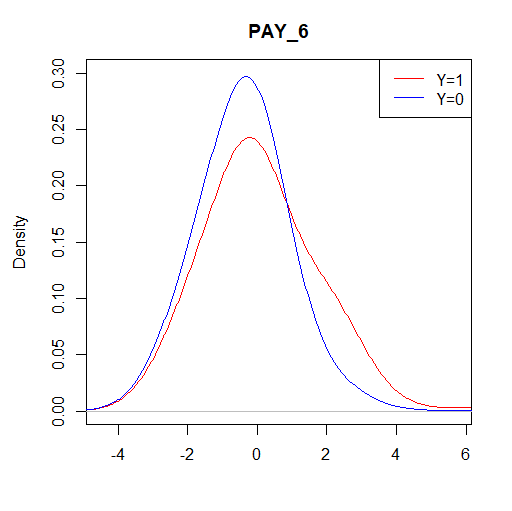} 
  
\vskip1.6em
\caption{UCI credit card data: We demonstrate the predictive performance and explainability of the LP-copula-based additive logistic regression model. The boxplots of the accuracy (AUC) are shown in the top left panel. On average, LP-logistic regression provides a 12\% boost in the accuracy. The LS-plot is shown in the top left. For easy interpretation, we display the scaled the LS-plot: $(\frac{\hal_{j1}}{max_j \hal_{j1}}, \frac{\hal_{j2}}{max_j \hal_{j2}})$. We can see a tight cluster around the origin $(0,0)$, which indicates that most of the variables are irrelevant for prediction (sparsity-assumption). Also, it is evident that the majority of the features have either location or scale information (differences), the only exception being the variable pay$\_$0---which denotes the repayment status of the last two months (-1=pay duly, 1=payment delay for one month, 2=payment delay for two months, and so on). This is further illustrated using three variables, as marked in the LS-plot; see also the two-sample density estimates shown at the bottom panel.}
\label{fig:uci}
\vspace{-.25em}
\end{figure}

\subsubsection{High-dimensional Copula-based Additives Logistic Regression} Generalize the univariate copula-logistic regression model \eqref{coplogthm} to the high-dimensional case
as follows:
\beq \label{eq:ALR}
\logit\big( \mu(x) \big)\,=\, \al_0 + \sum_{j=1}^p h_j(x_j).
\eeq
Nonparametrically approximate the unknown smooth $h_j$'s by LP-polynomial series of $X_j$
\beq \label{eq:ALR2}
h_j(x_j)\,=\,\sum_{k=1}^m \al_{jk} T_k(x_j;F_{X_j}), ~\text{for}\,j=1,\ldots,p.
\eeq

\begin{rem}[Estimation and Computation]
A sparse nonparametric LP-additive model of the form (\ref{eq:ALR}-\ref{eq:ALR2}) can be estimated by penalized regression techniques (lasso, elastic net, etc.), whose implementation is remarkably easy using \texttt{glmnet} R-function \citep{friedman2010regularization}:
$$\texttt{glmnet}(y \sim T_{\bf X}, 
\text{family=\texttt{binomial}})~~~~~~~$$
where $ T_{\bf X}$ simply is the column-wise stacked $[T_{X_1}\mid \cdots \mid T_{X_p}]$ LP-feature matrix. Another advantage of this formulation is that a large body of already existing theoretical work (see the monograph \cite{hastie2015book}) on  $\ell_1$-regularized logistic regression model can be directly used to study the properties of \eqref{eq:ALR2}.
\end{rem}

\vskip.1em
\begin{example}
\textit{UCI Credit Card data}. The dataset is available in the UCI Machine Learning Repository. It contains records of $n=30,000$ cardholders from an important Taiwan-based bank. For each customer, we have a response variable $Y$ denoting: default payment status (Yes = 1, No = 0), along with $p=23$ predictor variables (e.g., gender, education, age, history of past payment, etc.). We randomly partition the data into training and test sets, with an 80-20 split, repeated $100$ times. We measure the 
prediction accuracy using AUC (the area under the ROC curve). Fig. \ref{fig:uci} compares two kinds of lasso-logistic regressions: (i) usual version: based on feature matrix $X$; and (ii) LP-copula version: based on feature matrix $T_{{\bf X}}$.  As we can see,  LP-copula based additive logistic regression classifier significantly outperforms the classical logistic regression model. To gain further insight into the nature of impact of each variable, we plot the lasso-smoothed location and scale coefficients:
$$ \text{LS-Feature plot}:~~~~~\big(  \widehat{\al}_{j1},\, \widehat{\al}_{j2} \big),~~~~{\rm for}~j=1,\ldots,p.~~~~~~~~
$$
L-stands for location and S-stands for scale. The purpose of the LS-feature plot is to characterize `how' each feature impacts the classification task. For example, consider the three variables pay$\_$0, limit$\_$balance, and pay$\_$6, shown in the bottom panel Fig. \ref{fig:uci}.  Each one of them contains unique discriminatory information: pay$\_$0 has location as well as scale information, hence it appeared at the top-right of the LS-plot; The variable limit$\_$balance mainly shows location differences, whereas the variable pay$\_$6 shows contrasting scale in the two populations. In short, LS-plot explains `why and how' each variable is important using a compact diagram, which is easy to interpret by researchers and practitioners. 
\end{example}

\section{Conclusion: Copula-based Statistical Learning}
This paper makes the following contributions: 
(i) we introduce modern statistical theory and principles for maximum entropy copula density estimation that is self-adaptive for the \texttt{mixed}(X,Y)---described in Section \ref{sec:theory}. (ii) Our general copula-based formulation provides a unifying framework of data analysis from which one can systematically distill a number of fundamental statistical methods by revealing some completely unexpected connections between them. The importance of our theory in applied and theoretical statistics is highlighted in Section \ref{sec:Ut}, taking examples from different sub-fields of statistics: Log-linear analysis of categorical data, logratio biplot, smoothing large sparse contingency tables,  mutual information, smooth-$G^2$ statistic, feature selection, and  copula-based logistic regression. We hope that this new perspective on copula modeling will offer more effective ways of 
developing united statistical algorithms for mixed-($X,Y$).
\vskip.4em

\subsection*{Dedication: Two Legends from Two Different Cultures}
This paper is dedicated to the birth centenary of {\bf E. T. Jaynes} (1922--1998), the originator of the maximum entropy principle.  

I also like to dedicate this paper to the memory of {\bf Leo Goodman} (1928--2020)---a transformative legend of categorical data analysis, who passed away on December 22, 2020, at the age of 92 due to COVID-19.

This paper is inspired in part by the author's intention to demonstrate how these two modeling philosophies can be connected and united in some ways. This is achieved by employing a new nonparametric representation theory of generalized copula density. 

\vskip1em
\bibliographystyle{asa}
\bibliography{ref-bib}

\newpage
\section{Supplementary Appendix}
\renewcommand{\theequation}{6.\arabic{equation}}
\appendix
\renewcommand{\thesubsection}{A.\arabic{subsection}}
\subsection{Nonparametric LP-Polynomials}
\label{app:LPb}
\vskip.5em
{\bf Preliminaries}. For a random variable $X$ with the associated probability distribution $F_X$,  define the mid-distribution function as
$\Fm(x;F_X)=F_X(x)-\frac{1}{2}p(x;F_X)$ where $p(x; F_X)$ is probability mass function. The $\Fm(X;F_X)$ has mean $\Ex[\Fm(X;F_X)]=.5$ and  variance $\Var[\Fm(X;F_X)]=\frac{1}{12}\big( 1- \sum_x p^3(x;F_X) \big)$. Define first-order basis function $T_1(X;F_X)$ by standardizing $\Fm$-transformed random variable: 
\beq \label{eq:LP1st}
T_1(x;F_X)~=~\dfrac{\sqrt{12}\big\{\Fm(x;F_X) - 1/2\big\}}{\sqrt{1-\sum_x p^3(x;F_X)}}\eeq
where $\Ex[T_1(X;F_X)]=0.5$ and $\Var[T_1(X;F_X)]=1-\sum_x p^3(x)$.
Construct the higher-order LP-polynomial bases\footnote{Here the number of LP-basis functions $m$ is always less than $|\mathscr{U}|$, where $\mathscr{U}$ denotes the set of all unique values of $X$.} $\{T_j(X;F_X)\}_{1 \ge j\le m}$ by Gram-Schmidt orthonormalization of $T_1^2,T_1^3,\ldots,T_1^m$. We call these specially-designed polynomials of mid-distribution transforms as LP-basis, which by construction, satisfy the following orthonormality with respect to the measure $F_X$:  
\beq 
\int_x T_j(x;F_X) \dd F_X(x)=0,~~\text{and}~~\int_x T_j(x;F_X) T_k(x;F_X)\dd F_X(x)=\delta_{jk},
\eeq 
where $\delta_{jk}$ is the Kronecker delta.  Two particular LP-family of polynomials are given below:
\vskip.25em
~~~$\bullet$~{\bf gLP-Polynomial Basis}. Let $X \sim G$, where $G$ is a known distribution. Following the above recipe, construct parametric gLP-basis $\{T_j(X;G)\}_{j\ge 1}$ for the given distribution $G$. For $X \sim {\rm Bernoulli}(p)$, one can show that $T_1(x;G)=\frac{x-p}{\sqrt{p(1-p)}},$ $x=0,1;$ see \cite{D20copula} for more details. 
\vskip.25em

~~~$\bullet$~{\bf eLP-Polynomial Basis}. In practice, we only have access to a random sample $X_1,\ldots,X_n$ from an \textit{unknown} distribution $F$. To perform statistical data analysis, it thus becomes necessary to nonparametrically `learn' an appropriate basis that is orthogonal with respect to the (discrete) empirical measure $\wtF_X$. To address that need, construct LP-basis with respect to the empirical measure $\{T_j(X;\wtF_X)\}_{j\ge 1}$.
\vskip.25em
{\bf LP-unit Basis}. We will occasionally express the $T_j$'s in the quantile domain
\beq 
S_j(u;X)\,=\,T_j\big(  Q_X(u); F_X \big), ~0<u<1.
\eeq
We call these S-functions the unit LP-bases. The `S-form' and the `T-form' will be used interchangeably throughout the paper, depending on the context.

\vskip.25em
{\bf LP-product-bases and Copula Approximation}. The LP-product-bases $\big\{T_j(x;F_X)T_k(y;F_Y)\big\}$ can be used to expand any square-integrable function of the form $\Psi(F_X(x), F_Y(y))$. In particular, the logarithm of Hoeffding's dependence function \eqref{eq:depf} can be approximated by LP-Fourier series:
\beq \label{eq:sadepf}
\log \cop_{X,Y}\big(F_X(x),F_Y(y)\big) ~\sim~\sum_{j,k}\te_{j,k}T_j(x;F_X)T_k(y;F_Y).
\eeq
Represent \eqref{eq:sadepf} in the quantile domain by substituting $x=Q_X(u)$ and $y=Q_Y(v)$ to get the copula density expression \eqref{exp:model}.
\vskip.2em
\subsection{Two Cultures of Maximum Entropy Modeling}
\label{app:maxent}

The form of the maximum-entropy exponential model directly depends on the \textit{form} of the set of constraints, i.e., the sufficient statistics functions. The maxent distribution depends on the data only through the sample averages for these functions.

$\bullet$ Parametric maxent modeling culture: The traditional practice of maxent density modeling assumes that the appropriate sufficient statistics functions are known or given beforehand, which, in turn, puts restrictions on the possible `shape' of the probability distribution. This parametric maxent modeling culture was first established by Ludwig Boltzmann in 1877, and then later popularized by E. T. Jaynes in 1960s.

$\bullet$ Nonparametric maxent modeling culture: It proceeds by identifying a small set of most important sufficient statistics functions from data \citep{D11a2}. In the next step, we build the maxent probability distribution that agrees with these specially-designed relevant constraints. We have used LP-orthogonal polynomials to systematically and robustly design the constraining functions. See, Section A.3 for more discussion.

\subsection{Connection With Vladimir Vapnik's Statistical Invariants}
\label{app:VV}
Vladimir Vapnik (\citeyear{vapnik2020}) calls the sufficient statistics functions as `predicates' and the associated moment constraints as `statistical invariants.' While describing his learning theory he acknowledged that
\begin{quote}
``\textit{The only remaining question in the complete statistical learning theory is how to choose a (small) set of predicates. The choice of predicate functions reflects the intellectual part of the learning problem.}''
\end{quote}
The question of how to systematically design and search for informative predicates was previously raised by \cite{D11a2} in the context of learning maxent probability models from data; also see Section 2.2 of the main article where we discussed some concrete strategies to address this problem. However, it is important to note that \cite{vapnik2020} works entirely within the traditional least-square ($L_2$ risk) setup, instead of maximum-entropy framework.

\begin{figure}[t]
  \centering
\includegraphics[width=.7\linewidth,keepaspectratio,trim=4cm 1.5cm 1cm 1cm]{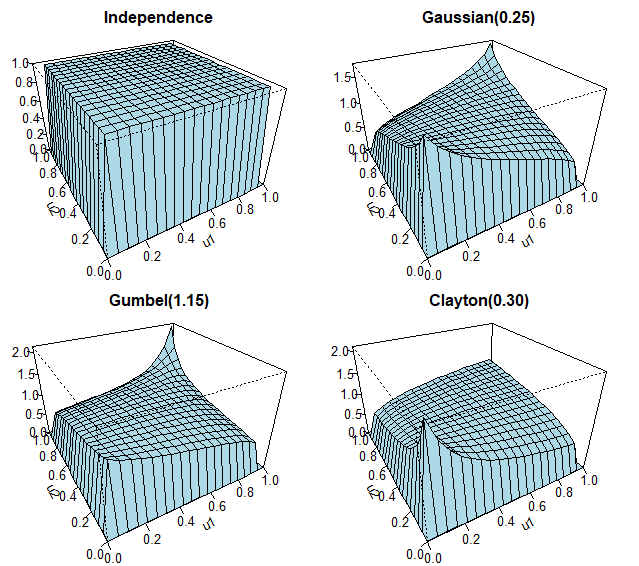}
\vskip4em
\caption{The copula models used for simulation study in Table \ref{tab:simulation}.}
\label{fig:copsimu}
\vspace{1.2em}
\end{figure}
\subsection{Empirical Power Study}
\label{app:simulation}
The simulation study is constructed as follows:  
\begin{table}[htb]
\vspace{2em}
\renewcommand{\arraystretch}{1.5}
    \centering
    \begin{tabularx}{.965\linewidth}{X@{\hskip 1.15in}X@{\hskip .7in}XXXXX}
\toprule
\multirow{2}*{Copula} & \multirow{2}*{Method} & \multicolumn{5}{c}{Values of $I$} \\   \cline{3-7}  
  &   &    5 & 20 & 40  & 50 & 100  \\
\midrule
\multirow{2}{*}{$U[0,1]^2$} & $~G^2$  & 0.058 &  \fbox{0.760}&  ~0 & ~0 & ~0\\
 & $~\widehat{G}^2$  &   0.062  & 0.044 & 0.059 &0.041&0.060\\
 \midrule
\multirow{2}{*}{Gaussian$(0.25)$} & $~G^2$  & 0.950 & 0.970&  ~0 & ~0 & ~0\\
 & $~\widehat{G}^2$& 0.920& 0.970 & 0.960 &0.940  &0.960\\
  \midrule
\multirow{2}{*}{Gumbel$(1.15)$} & $~G^2$  &0.790  & 0.960 &  ~0 & ~0 & ~0\\
 & $~\widehat{G}^2$  & 0.758& 0.930 & 0.910  &0.940  &0.90\\
   \midrule
\multirow{2}{*}{Clayton$(0.30)$} & $~G^2$  &  0.860& 0.950 &  ~0 & ~0 & ~0\\
 & $~\widehat{G}^2$    & 0.840 & 0.910 & 0.950  & 0.930  & 0.940\\
\bottomrule
\end{tabularx}
\vspace{.5em}

\caption{Empirical power study (of Sec \ref{sec:G2}) under four type of copula-dependency structure with $I\times I$ contingency tables; see Fig. \ref{fig:copsimu}. The LP-smoothed $G^2$ statistic is denoted by $\widehat{G}^2$. All tests were performed at nominal $\alpha=0.05$. Wilks’ $~G^2$ statistic starts to break down (both in terms of type-I error and power) even for mildly large tables. On the other hand, the proposed LP-smoothed  $\widehat{G}^2$-test continuous to work for contingency tables of all sizes ---from small-and-dense to large-and-sparse cases. It can thus be used as a ``default’’ nonparametric test of independence for categorical data.
} \label{tab:simulation}
\end{table}

Step 1. \textit{Constructing tables with different dependence structures}. We simulate $n=500$ independent random samples from the following copula distributions:
\begin{itemize}[itemsep=0pt,topsep=0pt]
\item Independent copula: $U[0,1]^2$;
\item Gaussian copula with  $\rho=0.25$;
\item Gumbel copula with parameter $\theta=1.15$.
\item Clayton copula with parameter $\theta=0.30$.
\end{itemize}
The shapes of the copulas are shown in Fig. \ref{fig:copsimu}. We convert each bivariate dataset into
$I\times I$ contingency table (having equally spaced bins) of counts $\{n_{ij}\}$, with $n=\sum_{i,j}n_{ij}$, $n_{i\bcdot} = \sum_j n_{ij}$, and $n_{\bcdot j} = \sum_i n_{ij}$.

Step 2. \textit{Simulating tables under null}. Generate $B$ (null) tables with given row marginal $\{n_{i\bcdot}\}$ and column marginal $\{n_{\bcdot j}\}$ using \citeauthor{patefield1981algR}'s (1981) algorithm. It is implemented in the \texttt{R}-function \texttt{r2dtable}.

Step 3. \textit{Power approximation}. We used $B=250$ null 
tables to estimate the 95\% rejection cutoffs at the significance level $\alpha=0.05$. The power is estimated based on $200$ independence tests for different sizes of contingency tables: $I=\{ 5, 20, 40, 50, 100\}$. We have fixed $m=4$ to compute the smooth-$G^2$ statistic, following Eqs. \eqref{eq:MIest} and \eqref{eq:sgMI}.
\vskip.5em

{\bf Result}. The type I error rates and power performances are reported in Table \ref{tab:simulation}. The $G^2$ test shows some strange Type-I error rate pattern: for medium-large contingency tables ($20\times 20$ case), it acts as an ultra-liberal test with Type I error ($0.760$, in box) much higher than the nominal level; while, for large-sparse cases (with $I \ge 40$) it behaves as an ultra-conservative test that yields Type I error rate much lower (almost zero) than the nominal level. On the other hand, $\widehat{G}^2$-test
maintains the type-I error rates close to the significance level $0.05$, even for highly sparse scenarios. In terms of power, the $G^2$ test is only reliable for the $5\times 5$ case. The higher power for the $20\times 20$ table is just a consequence of large type-I error phenomena. Conclusion: the conventional Wilks’ $~G^2$ is not a trustworthy test, even for moderately large tables, and should not be used blindly as a ``default’’ method. The smoothed $\widehat{G}^2$ performs remarkably well under all of these different scenarios and emerges as the clear winner.

\subsection{Additional Figures and Tables} \label{app:fig}

\begin{figure}[ h]
  \centering
\includegraphics[width=.47\linewidth,keepaspectratio,trim=4cm 1.5cm 1cm 1cm]{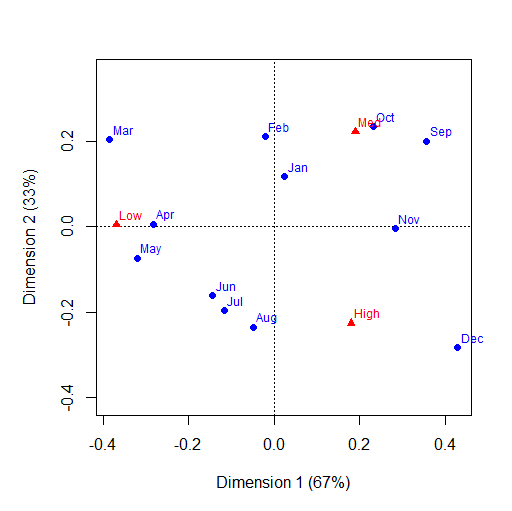}
\vskip2em
\caption{The classical correspondence analysis plot for the 1970 draft lottery data; implemented using the R-package \texttt{ca}. Contrast it with Fig. \ref{fig:draftUSA}(b).}
\label{fig:CA}
\end{figure}

\begin{figure}[ ]
\vskip2.4em
  \centering
\includegraphics[width=.46\linewidth,keepaspectratio,trim=3cm 2cm 2cm 4.5cm]{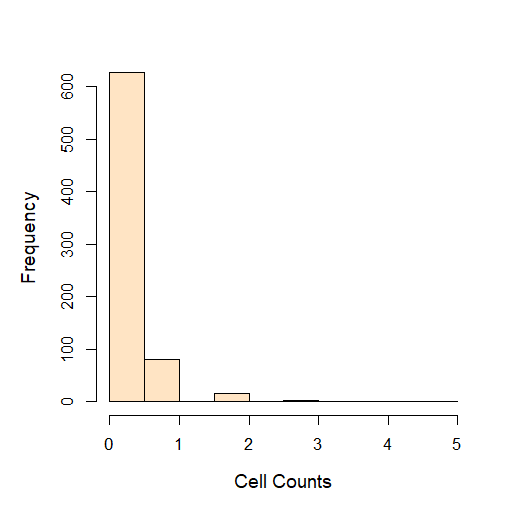}
\vskip2.5em
\caption{Zelterman data is a $28\times 26$ table, which is extremely sparse. Here we display the histogram of the observed cell frequencies. For more details see Example \ref{ex:zel}.}
\label{fig:appzelterman2}
\vspace{-.25em}
\end{figure}

\begin{figure}[ ]
  \centering
\includegraphics[width=.55\linewidth,keepaspectratio,trim=3cm 1cm 2cm 2cm]{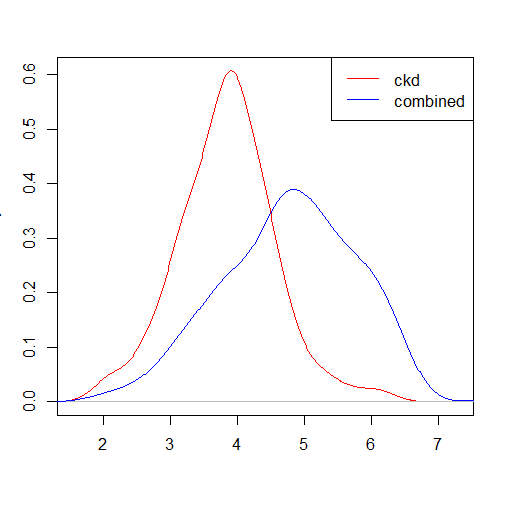}\\[2em]
\includegraphics[width=.55\linewidth,keepaspectratio,trim=3cm 2cm 2cm 2cm]{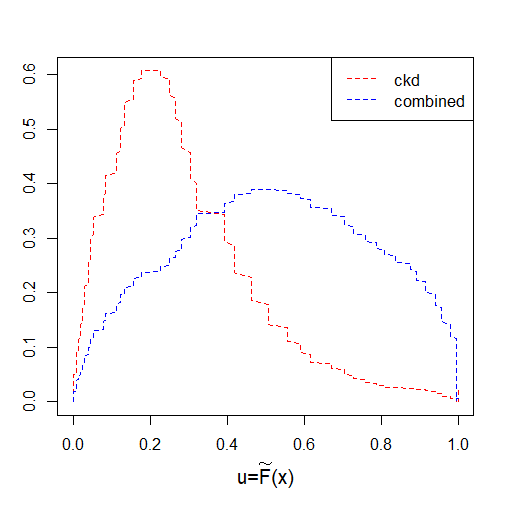}
\vskip3em
\caption{CKD data.}
\label{fig:ckd2}
\vspace{-.25em}
\end{figure}

\begin{figure}[ ]
  \centering
\includegraphics[width=.577\linewidth,keepaspectratio,trim=4cm 2cm 2.5cm 3cm]{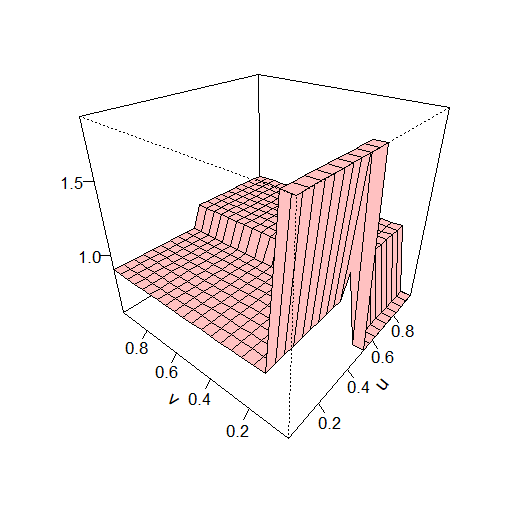}
\vskip.5em
\caption{LP-copula density estimate for Hellman $2\times 2$ data.}
\label{fig:hellmapp}
\end{figure}


\begin{table}[ht]
\setlength{\tabcolsep}{22pt}
\renewcommand{\arraystretch}{.82}
\centering
\begin{tabular}{cccc}
  \hline
 Months & High & Med & Low \\ 
  \hline
Jan &   9 &  12 &  10 \\ 
  Feb &   7 &  12 &  10 \\ 
  Mar &   5 &  10 &  16 \\ 
  Apr &   8 &   8 &  14 \\ 
  May &   9 &   7 &  15 \\ 
  Jun &  11 &   7 &  12 \\ 
  Jul &  12 &   7 &  12 \\ 
  Aug &  13 &   7 &  11 \\ 
  Sep &  10 &  15 &   5 \\ 
  Oct &   9 &  15 &   7 \\ 
  Nov &  12 &  12 &   6 \\ 
  Dec &  17 &  10 &   4 \\ 
   \hline
\end{tabular}
\caption{1970 US draft lottery data of Example \ref{ex:draft}. 
The lottery was conducted as follows:  $366$ possible birthdates (including Feb 29) were numbered as 1 through 366 on slips of paper; they
were put in a glass jar (after mixing them well) and drawn at random; the order in which they were selected was their drawing number or order of induction. Drawing numbers with 1-122 are marked as `high' risk category of being inducted; 123-244 as `medium' and 245-366 as `low' risk category.} \label{tab:draft}
\end{table}

\begin{table}[ht]
\centering
\begin{tabular}{rrrrrrrr}
  \hline
 & 0 & 1 & 2 & 3 & 4 & 5 & 6 \\ 
  \hline
0 & 0.17 & 0.15 & 0.07 & 0.02 & 0.01 & 0.00 & 0.00 \\ 
  1 & 0.11 & 0.12 & 0.08 & 0.01 & 0.03 & 0.01 & 0.00 \\ 
  2 & 0.03 & 0.04 & 0.03 & 0.02 & 0.01 & 0.00 & 0.01 \\ 
  3 & 0.02 & 0.01 & 0.02 & 0.02 & 0.00 & 0.01 & 0.00 \\ 
  4 & 0.00 & 0.00 & 0.01 & 0.01 & 0.00 & 0.00 & 0.00 \\ 
  7 & 0.00 & 0.01 & 0.00 & 0.00 & 0.00 & 0.00 & 0.00 \\ 
   \hline
\end{tabular}
\caption{Accident data. The raw empirical probability estimates $\widetilde{p}(x,y)$.} \label{tab:accEMP}
\end{table}

\begin{table}[ht]
\centering
\begin{tabular}{rrrrrrrr}
  \hline
 & 0 & 1 & 2 & 3 & 4 & 5 & 6 \\ 
  \hline
0 & 0.19 & 0.13 & 0.06 & 0.02 & 0.01 & 0.00 & 0.00 \\ 
  1 & 0.09 & 0.10 & 0.07 & 0.03 & 0.01 & 0.01 & 0.00 \\ 
  2 & 0.03 & 0.05 & 0.04 & 0.02 & 0.01 & 0.00 & 0.00 \\ 
  3 & 0.01 & 0.02 & 0.02 & 0.01 & 0.01 & 0.00 & 0.00 \\ 
  4 & 0.00 & 0.01 & 0.01 & 0.00 & 0.00 & 0.00 & 0.00 \\ 
  7 & 0.00 & 0.00 & 0.00 & 0.00 & 0.00 & 0.00 & 0.00 \\ 
   \hline
\end{tabular}
\caption{Accident data. The LP-smoothed probability estimates $\widehat{p}(x,y)$; see Example \ref{ex:shun}.} \label{tab:accLP}
\end{table}

\newpage

\end{document}